\documentclass[11pt]{amsart}                    

\usepackage{amssymb}  
\usepackage{graphicx}
\usepackage{setspace}
\input xy

\xyoption{all}

\theoremstyle{plain}

\newtheorem{theorem}{Theorem}

\newtheorem{lemma}{Lemma}[section]

\newtheorem{corollary}[lemma]{Corollary}

\newtheorem{fact}[lemma]{Fact}

\newtheorem{claim}[lemma]{Claim}

\newtheorem{question}[lemma]{Question}

\theoremstyle{definition}

\newtheorem{definition}[lemma]{Definition}

\newtheorem{remark}[lemma]{Remark}

\newtheorem{assumption}[lemma]{Assumption}

\newtheorem{notation}[lemma]{Notation}

\theoremstyle{remark}

\newcommand{\St}{\operatorname{\mathfrak{S}t}}

\newcommand{\cf}{\operatorname{cf}}

\newcommand{\LS}{\operatorname{LS}}
\newcommand{\lan}{\operatorname{L}}

\newcommand{\Aut}{\operatorname{Aut}}

\def\rg{\operatorname{rg}}

\newcommand{\ftp}{\operatorname{tp}}

\newcommand{\tp}{\operatorname{ga-tp}}

\newcommand{\gaS}{\operatorname{ga-S}}

\newcommand{\conc}{\Hat{\ }}

\newcommand{\aut}[1]{\operatorname{Aut}_{#1}(\mathfrak{C})}

\renewcommand{\phi}{\varphi}

\newcommand{\Union}{\bigcup}

\newcommand{\initial}\lessdot

\newcommand{\K}{\operatorname{\mathcal{K}}}

\newcommand{\id}{\operatorname{id}}
\newcommand{\s}{\operatorname{succ}}

\def\l{\langle}
\def\r{\rangle}

\newcommand{\C}{\mathfrak C}

\def\?{?\vadjust

{\vbox to 0pt{\vskip-7pt\hbox to 1.1\hsize{\hfill\huge ?!}}}}

\begin{document}

\title[Uniqueness of Limit Models in Classes with
Amalgamation]{Uniqueness of Limit Models in Classes with Amalgamation}

\author{Rami Grossberg}
\email[Rami Grossberg]{rami@cmu.edu}
\address{Department of Mathematics\\
Carnegie Mellon University\\
Pittsburgh PA 15213}

\author{Monica VanDieren}
\email[Monica VanDieren]{vandieren@rmu.edu}
\address{Department of Mathematics\\
Robert Morris University \\
Moon Township PA 15108}

 \author{Andr\'es Villaveces}
\email[Andr\'es Villaveces]{avillavecesn@unal.edu.co}
\address{Departamento de Matem\'aticas\\
Universidad Nacional de Colombia\\
Bogot\'a, Colombia}

 \thanks{
  The second author was partially sponsored for this work by grant
  DMS 0801313  of the \emph{National Science Foundation.}\\
  The third author was partially sponsored for this work by grant
  No. 1101-05-13605 of the \emph{Instituto Colombiano para el
    Desarrollo de la Ciencia y la Tecnolog\'\i a}, Colciencias.}

\date{January 15, 2015}

\maketitle

\begin{abstract}
 
  We prove:\\
{\bf Main Theorem:}  {\em Let $\K$ be an abstract
   elementary class satisfying the joint embedding and the amalgamation properties.  
Let $\mu$ be a cardinal above the the L\"owenheim-Skolem number of the class. 
Suppose $\K$ satisfies the disjoint amalgamation property for limit models
of cardinality $\mu$.
If $\K$   is $\mu$-Galois-stable, has no $\mu$-Vaughtian Pairs, does not have long splitting chains,
and satisfies locality of splitting,
then any two
$(\mu,\sigma_\ell)$-limits over
$M$ for  ($\ell\in\{1,2\}$) are isomorphic over $M$.}


This theorem extends results of Shelah from \cite{Sh 394},
\cite{Sh576}, \cite{Sh 600}, Kolman and Shelah  in \cite{KoSh} and
Shelah and Villaveces from \cite{ShVi}.  A preliminary version of our uniqueness theorem was used by Grossberg and VanDieren to prove a case of Shelah's
categoricity conjecture for tame abstract elementary classes in
\cite{GrVa2}. 
\end{abstract}

\section{Introduction}

We work in the general context of abstract elementary classes (AECs) with the
amalgamation property (AP), the disjoint amalgamation property, the joint
embedding property (JEP), and Galois-stability at one fixed cardinality
$\mu$ above the L\"owenheim-Skolem number. We prove the uniqueness of
limit models under a unidimensionality-like assumption of no $\mu$-Vaughtian pairs and superstability-like assumptions  of the
$\mu$-splitting dependence relation.

The basic model theory of abstract elementary classes (definitions, the role of the AP and the JEP, the existence of a ``monster model'' $\frak C$,
Galois types and the foundational development of stability theory in that
context) can be checked in the monograph~\cite{Gr2} and the
books~\cite{Ba},~\cite{Sh i}. For the sake of completeness, we
include some of the fundamentals of this
context here.

In 1977, Shelah, building on the work of Fra\"\i ss\'e and J\'onsson, identified
a non-elementary context in which a model theoretic analysis could be
carried out.  Shelah began to study classes of models, together with a
partial ordering of the class, which exhibit many of the properties that
the
models of a first order theory have with respect to the elementary
submodel relation.  Such classes were named abstract elementary
classes. They include classes of models axiomatizable in $L_{\omega_1,\omega}(\mathbf
Q)$. Both classification theory and stability theory may be carried
out to some extent within these classes. One strong advantage is that
there are no a priori compactness assumptions.
We reproduce the definition here. 
\begin{definition}\label{def AEC}
Let $\K$ be a class of structures
all in the same similarity type $\lan(\K)$, and
let $\prec_{\K}$ be a partial order on $\K$.  The ordered pair
$\langle \mathcal{K},\prec_{\mathcal{K}}\rangle$ is an
\emph{abstract elementary class,  AEC for short}
 iff
\begin{enumerate}

\item [A0] (Closure under isomorphism)
\begin{enumerate}

\item
For every $M\in \mathcal{K}$
 and every $\lan(\K)$-structure $N$ if $M\cong N$ then $N\in \mathcal{K}$.
\item Let $N_1,N_2\in\K$ and  $M_1,M_2\in \K$ such that
there exist $f_l:N_l\cong M_l$ (for $l=1,2$) satisfying
$f_1\subseteq f_2$ then  $N_1\prec_{\K}N_2$
implies that  $M_1\prec_{\K}M_2$.

\end{enumerate}

\item [A1]
For all $M,N\in \mathcal{K}$ if $M\prec_{\mathcal{K}} N$ then $M\subseteq N$.
 
\item [A2]
Let $M,N,M^*$ be $\lan(\mathcal{K})$-structures in $\K$.
If $M\subseteq N$, $M\prec_{\mathcal{K}} M^*$ and $N\prec_{\mathcal{K}}
M^*$, then
 $M\prec_{\mathcal{K}} N$.

\item [A3]
(Downward L\"owenheim-Skolem)
$\LS(\mathcal{K})$ is the minimal infinite
cardinal $\geq |\lan(\K)|$ 
such that
for every 
$M\in \mathcal{K}$
 and for every $A\subseteq |M|$ there exists $N\in \mathcal{K} $ such that
$N\prec_{\mathcal{K}} M, \; |N|\supseteq A$ and
$\|N\|\leq |A|+\LS(\mathcal{K})$.

\item [A4]
(Tarski-Vaught Chain)
\begin{enumerate}

\item
For every regular cardinal  $\mu$ and every 
$N\in \mathcal{K}$ if
$\langle M_i\in\K\mid M_i\prec_{\mathcal{K}} N\;,\;i<\mu\rangle$
is $\prec_{\K}$-increasing (i.e. $i<j\Longrightarrow M_i\prec_{\mathcal{K}}
M_j$) then
$\Union_{i<\mu}M_i\in \mathcal{K}$ and $\Union_{i<\mu}M_i\prec_{\mathcal{K}} N
$.

\item
For every regular $\mu$,
if  $\langle M_i\in\K_\mu\mid i<\mu\rangle$ is
$\prec_{\K}$-increasing then $\Union_{i<\mu}M_i\in \mathcal{K}$ and
$M_0\prec_{\mathcal{K}} \Union_{i<\mu}M_i$.
\end{enumerate}
\end{enumerate}

For $M$ and $N\in \K$ a monomorphism $f:M\to N$ is called an
\emph{$\K$-embedding} iff $f[M]\prec_{\K}N$.  Thus, $M\prec_{\K}N$ is
equivalent to ``$\id_M$ is a $\K$-embedding from $M$ into $N$.''

For $M_0\prec_{\K}M_1$ and $N\in\K$, the formula
$f:M_1\underset{M_0}{\to}N$ stands for $f$ is a $\K$-embedding
such that $f\restriction M_0=\id_{M_0}$.
\end{definition}

For a class $\K$ and a cardinal $\mu\geq\LS(\K)$ let
\[
\K_\mu:=\{M\in
\K \;:\; \|M\|=\mu\}.
\]

In practice, abstract elementary classes were not as approachable as one
would hope and much work in non-elementary model theory takes place in
contexts which additionally satisfy the amalgamation property:

\begin{definition}\label{ap defn}
Let $\mu\geq \LS(\K)$.  We say that $\K$ has the \emph{$\mu$-amalgamation
property} (\emph{$\mu$-AP}) iff for any $M_\ell\in \K_\mu$ (for $\ell\in
\{0,1,2\}$) such that $M_0\prec_{\K} M_1$ and $M_0\prec_{\K}M_2$ there are
$N\in\K_\mu$ and $\K$-embeddings $f_\ell:M_\ell
\to N$ such that $f_\ell\restriction M_0=\id_{M_0}$ 
for $\ell =1,2$.


We say that $\K$ has the \emph{amalgamation property} (\emph{AP}) iff any
triple of models from $\K_{\geq\LS(\K)}$ can be amalgamated.
\end{definition}

\begin{remark}
\begin{enumerate}

\item
Using the isomorphism axioms we can see that $\K$ has the
$\lambda$-AP iff for any $M_\ell\in\K_\lambda$ (for  $\ell\in\{0,1,2\}$)
such that
$M_0\prec_{\K}M_\ell$ (for $\ell\in\{1,2\}$) there are $N\in\K_\lambda$
and $f:M_1\underset{M_0}{\to}N$ such that $N\succ_{\K}M_2$.

\item
Using the axioms of AECs it is not difficult to prove that if $\K$ has
the $\lambda$-AP for every $\lambda\geq\LS(\K)$ then $\K$ has the AP.

\end{enumerate}
\end{remark}

The roots of the following fact can be traced back to J\'onsson's 1960
paper~\cite{Jo}; the present formulation is from \cite{Gr1}:

\begin{fact}
Let $\l\K,\prec_{\K}\r$ be an AEC with no maximal models and suppose that there is $\lambda\geq\kappa >\LS(\K)$ such
that $K_{<\lambda}$ has the AP and the
JEP. Suppose $M\in\K$.
If  
 $\lambda^{<\kappa}=\lambda\geq\|M\|$ then there
exists 
$N\succ M$ of cardinality $\lambda$ which
is $\kappa $-model-homogeneous.
\end{fact}

Thus if an AEC $\K$ has AP and JEP, then
like in first-order stability theory we may assume that 
there is a large  model-homogeneous $\C\in \K$ that acts like a monster
model. 

We will refer to the model $\mathfrak C$ as the {\it monster model}.
All models considered will be of size less than $\|\mathfrak C\|$, and we
will find realizations of types we construct inside this monster
model.
From now on, we assume that the monster model $\mathfrak C$ has been fixed.  We use the notation $\aut{M}$ to denote the set of automorphisms of $\C$ fixing $M$ pointwise.  

%
%
%
%

The notion of type as a set of formulas, even when the class is described
in some infinitary logic, does not 
behave as nicely as in first-order logic.
A replacement was introduced by Shelah in \cite{Sh
394}.  In order to avoid confusion between this  and the classical, syntactic notion, we will use the terminology  in
\cite{Gr2} and call this alternative notion the
\emph{Galois type}.

Since in this paper we deal only with AECs with the AP property,
the notion of Galois type has a simpler definition than in the general
case.

\begin{definition}[Galois types]\label{type.def}
Suppose that $\K$ has the AP. \hfill
\begin{enumerate}
\item
Given $M\in\K$ consider the action of $\aut M$ on $\C$, for an element
$a\in |\C|$ let $\tp(a/M)$ denote the \emph{Galois type of $a$ over
$M$} which is defined as the orbit of $a$ under $\aut M$.

\item
For $M \in \K$, we let
\[
\gaS(M) = \{\  \tp(a/M) \;:\; a\in |\C|
  \}.
\]

\item  $\K$ is \emph{$\lambda$-Galois-stable}
\index{$\lambda$-Galois-stable}\index{Galois-stable} iff
\[
N\in\K_\lambda\implies |\gaS(N)|\leq\lambda.
\]

\item
Given $p \in \gaS(M)$ and $N \in \K$ such that $N\succ_{\K}M$, we say that
$p$ is \emph{realized} by $a \in N$ iff $\tp(a/M)= p$.
Just as in the first-order case we will write $a\models p$ when $a$ is a
realization of $p$.
\item For $h\in\Aut(\C)$ and $p=\tp(a/M)$, then the notation $h(p)$ refers to $\tp(h(a)/h(M))$.
\end{enumerate}
\end{definition}

For a more detailed discussion of Galois types, their extensions,
restrictions, equivalent forms and generalizations, the reader may
consult~\cite{Gr2}.

While the amalgamation property is useful for dealing with Galois types,  in this paper we require a stronger version of AP for one of the steps in the proof of the uniqueness of limit models.  Specifically, we use $\mu$-disjoint amalgamation over limit models to prove that relatively full towers are limit models (Theorem \ref{relatively full is limit}).  

\begin{definition}\label{def:DAP}  Let $\K$ be an abstract elementary class.
$\K$ has the $\mu$-\emph{disjoint amalgamation  property}  
\index{Disjoint Amalgamation  Property}
\index{$\lambda$-Disjoint Amalgamation  Property} ($\mu$-DAP) iff for every $M_\ell\in \K_\mu$ (for 
$\ell =0,1,2$) such that 
$M_0\prec_{\K} M_\ell$ (for $\ell=1,2$) 
there are $N\in \K_\mu$ which is a 
$\K$-extension of $M_2$
and
a $\K$-embedding $f:M_1\underset{M_0}{\to} N$ 
such that 
$ f[M_1]\cap M_2=M_0$.

We say that a class has the \emph{disjoint amalgamation property} iff it has the
 $\mu$-disjoint amalgamation property 
for every $\mu\geq \LS(\K)$.  We write DAP for short.  In this paper we only require that disjoint amalgamation hold for the subclass of all limit models of $\K_\mu$.
\end{definition}



%


The next notion to consider is that of a saturated model.  
In homogeneous abstract
elementary classes (see, for example, \cite{GrLe}) where one may study
classes of models omitting given sets of types, the existence of
a saturated model presents some problems.  One solution is to consider models which realize as many types as possible.  Such models are called Galois-saturated.  More formally, a model $M$ of size $\kappa >\LS(\K)$ is \emph{Galois-saturated} if it realizes all Galois types
over submodels $N\prec_{\K} M$ of  cardinality $<\kappa$.
When
stability theory has been ported to
contexts more general than first order logic, many situations have
appeared when Galois-saturated models do not fulfill the main roles that saturated models  play
in elementary classes.

The main concept of this paper is Shelah's {\it limit model} which
(among other things) serves as a substitute for the role of saturation in stability
theory (see~\cite{Gr2},\cite{ShVi},\cite{Sh i}, etc.) or at least serves as a stepping stone to prove the properties of Galois-saturated models.  
For example, 
under the assumption of
categoricity with reasonable stability conditions, the
existence of Galois-saturated models in singular cardinals is not straightforward and is
proved by first considering limit models \cite{Sh 394}.  In some contexts limit models have been successfully used
as ``tools'' towards finding Galois-saturated models (\cite{KoSh} and \cite{Sh 472}). Furthermore, the
notion of limit model refines the notion of saturation; more detailed
information is given on the particular way one model is embedded
inside another.

Limit models appear in  \cite{KoSh} and in \cite{Sh576}  under the name \emph{$(\mu,\alpha)$-saturated
  models}.    In \cite{Sh 600}, Shelah calls this notion \emph{brimmed}.
Later papers, beginning with Shelah-Villaveces \cite{ShVi}, adopt the name
\emph{limit models}.  We use the more recent terminology.
Before defining limit models, we must introduce their building blocks, universal extensions.


\begin{definition}\label{def:universal over}
\begin{enumerate}
\item \index{universal over!$\kappa$-universal over}
Let $\kappa$ be a cardinal $\geq \LS(\K)$.
We say $M^*\succ_{\K} N$ is \emph{$\kappa$-universal over $N$} iff 
for every $N'\in\K_{\kappa}$ with $N\prec_{\K}N'$ there exists
a $\K$-embedding $g:N'\underset{N}{\to} M^*$ such that the following diagram commutes:

\[
\xymatrix{\ar @{} [dr] N'
\ar@{.>}[dr]^{g}  & \\
N \ar[u]^{\id} \ar[r]_{\id} 
& M^*  
}
\]

\item \index{universal over}

We say $M^*$ is \emph{universal over $N$} or $M^*$ is \emph{a universal
extension of $N$} iff 
$M^*$ is $\|N\|$-universal over $N$.

\end{enumerate}
\end{definition}

\begin{definition}\label{def:limitmodels}[Limit models]
Consider $\mu\geq\LS(\K)$ and $\alpha<\mu^+$ a limit ordinal and $N\in\K_\mu$.
We say that $M$ is \emph{$(\mu,\alpha)$-limit model over $N$}
\index{$\alpha$-limit over $N$}
\index{limit over $N$}
iff there exists an increasing and continuous chain $\langle M_i\in\K_\mu\mid
i<\alpha\rangle$ such that 
$M_0=N$; $M=\bigcup_{i<\alpha}M_i$; $M_i$ is a proper $\K$-submodel of $M_{i+1}$;
and $M_{i+1}$ is universal over $M_i$
for all
$i<\alpha$.
\end{definition}

From Theorem \ref{exist univ} we get that for $\alpha\leq\mu^+$ there always
exists a $(\mu,\alpha)$-limit model  provided $\K$
has the AP, has no maximal models and is $\mu$-Galois-stable.  This theorem was stated without proof as Claim 1.16 in \cite{Sh 600},
for a proof see \cite{GrVa1} or \cite{Gr1}.

\begin{theorem}[Existence]\label{exist univ}
Let $\K$ be an AEC without maximal models and suppose it is
Galois-stable in
$\mu$.  If
$\K$ has the amalgamation property  then for every
$N\in\K_\mu$ there exists
$M^*\succneqq_{\K}N$, universal over $N$ of cardinality $\mu$.
\end{theorem}

The following theorem partially clarifies the analogy with saturated
models:

\begin{theorem}\label{limitfirstorder}
Let $T$ be a stable, complete, first-order theory and let $\K$ be the elementary
class of models of $T$ with the usual notion of elementary submodel.
If $M$ is a $(\mu,\delta)$-limit model for
$\delta$ a limit ordinal with $\cf (\delta)\geq \kappa(T)$, then $M$ is saturated. 

\end{theorem}

\begin{proof}
 Use an argument similar to the proof
 of~\cite[Theorem III 3.11]{Sh e}.
\end{proof}

Thus 
in elementary classes superstability implies that  limit
models are saturated, in particular are unique.
 This raises the following natural question for AECs:

\begin{question}[Uniqueness problem]\label{q:Uniqueness problem}
Let $\K$ be an AEC, $\mu\geq \LS(\K)$,  $M\in\K_\mu$ and $\sigma_1,\sigma_2$ limit ordinals $<\mu^+$,
 and suppose that for $\ell=1,2$, $N_\ell$ is a $(\mu,\sigma_\ell)$-limit model
over $M$. What ``reasonable'' assumptions on
$\K$ will imply that $\exists f:N_1\cong_M N_2$?
\end{question}

Question \ref{q:Uniqueness problem} is non-trivial only for the
case where
$\cf(\sigma_1)\neq\cf(\sigma_2)$.  Using a back and forth
argument one can show
that when $\cf (\sigma_1)=\cf( \sigma_2)$, we get uniqueness without any
assumptions on $\K$.   More
precisely:

\begin{fact}
\label{unique limits}
Let $\mu\geq \LS(\K)$ and $\sigma<\mu^+$.  
If $M_1$ and $M_2$ are $(\mu,\sigma)$-limits over $M$, then
there exists an isomorphism $g:M_1\to M_2$ such that
$g\restriction M = \id_M$.  Moreover if $M_1$ is a $(\mu,\sigma)$-limit
over $M_0$, if $N_1$ is a $(\mu,\sigma)$-limit over $N_0$ and
if $g:M_0\cong N_0$, then there exists a $\K$-embedding, $\hat g$,
extending $g$ such that $\hat g:M_1\cong N_1$.
\end{fact}


\begin{fact}
\label{sigma and cf(sigma) limits}\index{uniqueness of limit models!of the
same cofinality} Let $\mu$ be a cardinal and $\sigma$ a limit ordinal with
$\sigma<\mu^+
$.  If $M$ is a $(\mu,\sigma)$-limit
model, then $M$ is a $(\mu,\cf(\sigma))$-limit model.
\end{fact}

The  main result of this paper provides an answer to Question
\ref{q:Uniqueness problem}.

\begin{theorem}[Main Theorem]\label{main Theorem}  Let $\K$ be an AEC
  without maximal models, and
$\mu>\LS(\K)$.
Suppose $\K$ satisfies the $AP$ and $JEP$ and the subclass of limit
models of $\K$ satisfies $\mu$-DAP.  If $\K$   is
$\mu$-Galois-stable, does not have long splitting chains, has no
$\mu$-Vaughtian pairs and
satisfies locality of splitting\footnote{See Assumption \ref{splitting
assm} for the precise description of long splitting chains and locality.},
then any two
$(\mu,\sigma_\ell)$-limits over
$M$ for  ($\ell\in\{1,2\}$) are isomorphic over $M$.

\end{theorem}

\begin{remark}\label{Drueck}
After reading preprints of this paper, Fred Drueck in his Ph.D. thesis \cite{Dr} pointed out that the disjoint amalgamation property is not necessary to carry out the arguments here.  
In particular, it is not needed in Theorem \ref{relatively full is limit}.
We leave the assumption in this paper for historical accuracy.

\end{remark}

The last section of this paper (see pages~\pageref{context analysis}
and ff.)
describes different approaches to Question \ref{q:Uniqueness problem}.

We thank John Baldwin, Tapani Hyttinen and Pedro Zambrano for helping
to clarify the presentation. 
We also thank the referee for valuable suggestions, remarks and
an example of an $\aleph_1$-categorical AEC failing DAP over countable
models but having DAP over limit models.

\section{The Setting}

In what follows, $\K$ is assumed to be an AEC, and $\mu$ is a cardinal
$\geq \LS(\K)$.  In this section we summarize all of the assumptions that will be made on the class $\K$, and in the subsequent sections we introduce two of the main components of the proof of the uniqueness of limit models:  strong types and towers.

 We will prove the uniqueness of limit models in $\mu$-Galois stable AECs that are essentially unidimensional and are
 equipped with a moderately well-behaved dependence relation.  
  We will use $\mu$-splitting as the dependence relation, but  any
dependence relation which is local and has
 existence, uniqueness and extension properties suffices.
 
\begin{definition}\index{$\mu$-splits}\index{Galois-type!$\mu$-splits} 
A
type
$p\in \gaS(M)$ \emph{$\mu$-splits} over $N\in \K_{\leq \mu}$ if and only if
$N$ is a $\prec_{\K}$-submodel of $M$ and
 there exist $N_1,N_2\in\K_{\mu}$ and a $\K$-mapping $h$  such  that
$N\prec_{\K}N_l\prec_{\K}M$ for $l=1,2$ and
$h:N_1\to N_2$ with $h\restriction N=\id_N$ and 
$p\restriction N_2\neq h(p\restriction N_1)$.
\end{definition}

The existence property for non-$\mu$-splitting types follows from Galois
stability in $\mu$:

\begin{fact}[Existence - Claim 3.3 of \cite{Sh 394}]\label{nonsplit thm}
  Assume $\K$ has AP and is
  Galois-stable in $\mu$. For every $M\in\K_{\geq\mu}$ and $p\in \gaS(M)$, there
  exists $N\in\K_\mu$ such that $p$ does not $\mu$-split over $N$.
\end{fact}

The uniqueness and extension properties of  non-$\mu$-splitting types hold
for types over limit models:
\begin{fact}[Uniqueness - Theorem I.4.15 of \cite{Va1}]\label{unique ext}
Let $N\prec_{\K} M\prec_{\K} M'$ be models in $\K_\mu$ such that
$M'$ is universal over
$M$ and
$M$ is universal over $N$.  If $p\in \gaS(M)$ does not $\mu$-split over
$N$, then there is a unique $p'\in\gaS(M')$ such that $p'$ extends $p$
and $p'$ does not $\mu$-split over $N$.
\end{fact}

A variation of this fact is later used in an induction construction in
the proof of Theorem \ref{reduced are cont}.  We state it explicitly
here:

\begin{fact}[Theorem I.4.10 of \cite{Va1}]\label{splitting extension lemma}
Let $M,N,M^*$ be models in $\K_\mu$.  Suppose that $M$ is universal
over $N$ and that $M^*$ is an extension of $M$.  If a type $p=\tp(a/M)$
does not $\mu$-split over $N$ then there exists an automorphism $g$ of
$\C$ fixing $M$ such that $\tp(g(a)/M^*)$ does not $\mu$-split over
$N$ and $\tp(g(a)/M)=p$.
\end{fact}

The other concepts that show up in the assumptions of the main theorem of this paper are minimal types and $\mu$-Vaughtian Pairs.
\begin{definition}
\begin{enumerate}
\item For $M$ a model of cardinality $\mu$, $p\in\gaS(M)$ is \emph{minimal} if it is non-algebraic and for each $N$ extending $M$ of cardinality $\mu$ there is a unique non-algebraic extension of $p$ to $N$.
\item
For $M$ a limit model of cardinality $\mu$ a \emph{$\mu$-Vaughtian Pair} is a pair of models $M'$ and $N'$ of cardinality $\mu$ if $M\preceq_{\K}M'\prec_{\K}N'$ and 
if there exists $p\in\gaS(M)$ a minimal type so that $N'$ contains no new realizations of $p$, in other words, $p(M')=p(N')$.

\end{enumerate}
\end{definition}

\begin{fact}[Existence of minimal types - reference \cite{Sh 394}]\label{existence of minimal}
Let $\mu>\LS(\K)$.  If $\K$ is Galois-stable in $\mu$, then for every $M\in\K_{\mu}$ and every $q\in\gaS(M)$, there are $N\in\K_{\mu}$ and $p\in\gaS(N)$ such that $M\preceq_{\K}N$, $q\leq p$ and $p$ is minimal.
\end{fact}

\begin{fact}[Claim $(*)_8 $ of Theorem 9.7 of \cite{Sh 394}]\label{no vp fact}
If $\K$ is categorical in some successor cardinal $\lambda^+>\LS(\K)^{+}$, then for every $\mu$ satisfying $\LS(\K)\leq\mu\leq\lambda$, there are no $\mu$-Vaughtian Pairs.
\end{fact}


It is worth mentioning that our ``no $\mu$-Vaughtian
pairs'' assumption is much weaker in general than assuming
categoricity (as in earlier version of the proof): even in First
Order, theories such as the theory of Real Closed Fields are quite far
from being categorical but also have no Vaughtian pairs. Of course,
under $\omega$-stability, no Vaughtian pairs and categoricity are
equivalent (in First Order). But our stability assumptions are of
``superstable'' nature - under these, categoricity is quite stronger
than no $\mu$-Vaughtian pairs.

Here are the assumptions of the paper:

\begin{assumption}\label{splitting assm}
  $\K$ is an AEC with  the $\mu$-DAP\footnote{See Remark \ref{Drueck} which indicates only the amalgamation property is necessary.} over limit models and JEP, and $\K$ satisfies the
  following properties:
\begin{enumerate}
\item\label{monster assm} All models are submodels of a fixed monster
  model $\frak C$.\footnote{Notice that this already implies the full AP.}
\item\label{stable assm} $\K$ is stable in $\mu$.
 \item\label{unidimensionality}  There are no $\mu$-Vaughtian Pairs.
\item\label{split assm} $\mu$-splitting in $\K$ satisfies the following
  locality (sometimes called continuity) and ``no long splitting chains''
  properties.\\
%
%
For all infinite $\alpha$, for every sequence $\langle M_i\mid i<\alpha\rangle$ of
  limit models of cardinality $\mu$ and for every $p\in\gaS(M_\alpha)$, where
  $M_\alpha=\bigcup_{i<\alpha}M_i$, we have that
\begin{enumerate}
\item\label{locality} If for every $i<\alpha$, the type $p\restriction
  M_i$ does not $\mu$-split over $M_0$, then $p$ does not $\mu$-split over
  $M_0$.
\item\label{no long splitting chain} There exists $i<\alpha$ such that $p$
does not $\mu$-split over $M_i$.
\end{enumerate}
\end{enumerate}
\end{assumption}

In the context of an AEC with the full amalgamation
  property and JEP, categoricity in a cardinal $\lambda > \mu$ implies all
  parts of Assumption \ref{splitting assm}. For a proof of Assumption
  \ref{splitting assm}.\ref{stable assm} from categoricity, see Claim
  1.7 of \cite{Sh 394} or \cite{Ba}. The observation that  assumption
  \ref{splitting assm}(\ref{locality}) follows from categoricity is a
  consequence of Observation 6.2 and Main Lemma 9.4 of \cite{Sh 394}.
  Lemma 6.3 of \cite{Sh 394} is the statement that assumption
  \ref{splitting assm}(\ref{no long splitting chain}) follows from
  categoricity when the cofinality of the categoricity cardinal is
  larger than $\mu$.

  Assumption \ref{splitting assm} also holds in contexts without the
  assumption of categoricity.  First let us consider $\mu$-DAP.  The 
  $\mu$-DAP over limit models holds for free in first order classes of
  the form $(Mod(T),\prec)$
  for complete $T$.  As the referee has pointed out, in AECs, $\mu$-DAP does not generally hold over arbitrary models.  Consider the class $\K$ of
    structures with two sorts $U$ and $V$ and a binary relation $<$ on
    $U$ such that for each model $M$, $U^M$ is well-ordered by $<^M$
    with order type at most $\omega$, $V^M$ is empty when $U^M$ is
    finite and if non-empty, $V^M$ is infinite.  By defining $\prec_{\K}$ by $M\prec_{\K}
    N$ iff  $U^M$ is an initial segment of $U^N$ and $V^M\subset V^N$, we get an
AEC with $LS(\K)$ equal to $\aleph_0$.  $\K$ satisfies the AP and JEP and is
$\aleph_1$-categorical.  It fails to have the $\aleph_0$-DAP yet has
$\aleph_0$-DAP over limit models.

However, there are AECs in which $\mu$-DAP does hold.
DAP 
holds in homogeneous classes (see
  \cite{Sh 3} or \cite{Po}), in excellent classes (see \cite{Sh87b})
  and is an axiom in the definition of finitary classes (see
  \cite{HyKe}). It also holds for cats consisting of existentially
  closed models of positive Robinson theories (\cite{Za}).  In each of
  these contexts dependence relations satisfying Assumption
  \ref{splitting assm} have been developed. Finally, the locality and
  existence of non-$\mu$-splitting extensions are akin to consequences
  of superstability in first order logic.

\begin{fact}[``No long splitting chains'' follows from stability in FO]
  Suppose that $T$ is first order complete. If $T$ is stable
  then Assumption~\ref{splitting assm}(\ref{no long splitting chain})
  holds for $\alpha$ such that $\cf(\alpha)\geq |T|^+$.
\end{fact}

\begin{proof}
  Let $\langle M_i|i\leq \alpha\rangle$ be an increasing sequence of saturated models and  $p\in
  S(M_\alpha)$ be such that $\forall i<\alpha$, $p$ $\mu$-splits over $M_i$. Let
  $\varphi_i(\bar{x},\bar{y})$ be a formula witnessing the splitting of
  $p\restriction M_{i+1}$ over $M_i$. As $\cf(\alpha)\geq |T|^+$,
  there exists $S\subset \alpha$
  infinite such that $i,j\in S\Rightarrow \varphi_i=\varphi_j$.

  Without loss of generality, suppose that $\langle M_n|n\leq
  \omega\rangle$ is an increasing sequence of
  saturated models, and $p\in S_\varphi(M_\omega)$ is such that
  $\bar{a}_i,\bar{b}_i\in M_{i+1}$ witness that $p\restriction M_{i+1}$
  splits over $M_i$. Then
  $p(x_1,\bar{y}_1,\bar{z}_1,x_2,\bar{y}_2,\bar{z}_2)$ and $\{
  \bar{d}_i|i<\omega\}$ witness that $p$ has the order property, where
  $\bar{d}_i=\bar{a}_i\conc \bar{b}_i\conc c_i$, $c_i\in M_{i+2}$ and
  \[ c_i\models p\restriction \{ \bar{a}_k,\bar{b}_k|k\leq i\}\cup \{
  d_k|k<i\}.\]
  Now use~\cite[Lemma VII, 2.12]{Gr1}.
\end{proof}

\section{Strong Types}
Under the assumption of $\mu$-stability, we can define
 \emph{strong types}
 as in \cite{ShVi}. These strong types will allow us to achieve a
 better control of extensions of towers of models
than what we obtain using just Galois types.

\begin{definition}[Definition 3.2.1 of \cite{ShVi}]\label{strong type
    defn}

For $M$ a $(\mu,\theta)$-limit model, \index{strong
    types}\index{Galois-type!strong}\index{$\St(M)$}\index{$(p,N)$}
    let
$$\St(M):=\left\{\begin{array}{ll}
(p,N)
& 
\left|\begin{array}{l}
N\prec_{\K}M;\\
N\text{ is a }(\mu,\theta)\text{-limit model};\\
M\text{ is universal over }N;\\
p\in \gaS(M)\text{ is non-algebraic}\\
\text{and }p\text{ does not }\mu\text{-split over }N.
\end{array}\right\}
\end{array}\right .
$$
Elements of $\St(M)$ are called {\em strong types.}
Two strong types $(p_1,N_1)\in\St(M_1)$ and $(p_2,N_2)\in\St(M_2)$
are
\emph{parallel} iff for every $M'$ of cardinality $\mu$ extending $M_1$
and $M_2$ there exists $q\in\gaS(M')$ such that $q$ extends both $p_1$
and $p_2$ and $q$ does not $\mu$-split over $N_1$ nor over $N_2$. 

\end{definition}

\begin{remark}
  Under the assumption of the existence of universal extensions, it is
  equivalent to say two strong types $(p_1,N_1)\in\St(M_1)$ and
  $(p_2,N_2)\in\St(M_2)$ are parallel iff for some $M'$ of
  cardinality $\mu$ universal over some common extension of $M_1$ and
  $M_2$ there exists $q\in\gaS(M')$ such that $q$ extends both $p_1$ and
  $p_2$ and $q$ does not $\mu$-split over $N_1$ and $N_2$.
\end{remark}

\begin{lemma}[Monotonicity of parallel types]\label{monotonicity of
    parallel}
  Suppose $M_0,M_1
\in\K_\mu$ and $M_0\prec_{\K}M_1$ and $(p,N)\in\St(M_1)$.  If
$M_0$ is universal over $N$, then we have
$(p\restriction M_0,N)$ is parallel to $(p,N)$.
\end{lemma}
\begin{proof}
Straightforward using the uniqueness of non-$\mu$-splitting extensions.
\end{proof}

\begin{notation}
Let $M,M'\in\K_\mu$ and
suppose that $M \prec_{\K} M'$.
For $(p,N)\in \St(M')$, if $M$ is universal over $N$, we define the
restriction
$(p,N)\restriction M\in \St(M)$\index{$(p,N)\restriction M$} to be
$(p\restriction M,N)$.
If we write $(p,N)\restriction M$, we mean that $p$ does not $\mu$-split
over $N$ and $M$ is universal over
$N$.
We denote by $\sim$ the parallelism relation between strong types
in $\St(M)$, for fixed $M$.
\end{notation}

Notice that $\sim$ is an equivalence relation on $\St(M)$ (see
\cite{Va1}).  Stability in $\mu$ implies that there are few strong types over any model of cardinality $\mu$:

\begin{fact}[Claim 3.2.2 (3) of \cite{ShVi}]\label{St small}
If $\K$ is Galois-stable in $\mu$, then for any $M\in \K_\mu$,
$|\St(M)/\sim|\leq\mu$.
\end{fact}

\section{Towers}
To each $(\mu,\theta)$-limit model $M$ we can naturally associate a
$\prec_{\K}$-increasing chain  $\bar M=\langle M_i\in\K_\mu\mid
i<\theta\rangle$ witnessing that $M$ is a
$(\mu,\theta)$-limit model (that is, $\Union_{i<\theta}M_i=M$ and $M_{i+1}$ is
universal over $M_i$).  Furthermore, by Facts \ref{unique limits} and
\ref{sigma and cf(sigma) limits}  we can require that this chain
satisfies additional requirements such as $M_{i+1}$ is a limit model
over $M_i$. In this section we will be considering a related
chain of models which we will refer to as a tower (see
Definition \ref{def:towers}). But first, we will describe how towers
will be used to prove the main theorem of this paper.

To prove the uniqueness of limit models we will construct a model
which is simultaneously a $(\mu,\theta_1)$-limit model over some
fixed model $M$ and a $(\mu,\theta_2)$-limit model over $M$.  Notice that,
by Fact \ref{unique limits}, it is enough to construct a model $M^*$ that
is simultaneously a $(\mu,\omega)$-limit model and a $(\mu,\theta)$-limit
model for arbitrary ordinal $\theta<\mu^+$. By Fact \ref{sigma and cf(sigma)
  limits} we may assume that $\theta$ is a limit ordinal $<\mu^+$ such that
$\theta=\mu\cdot\theta$.

So, we actually construct an array of models with $\omega +1$ rows and the
number of columns of this array will have the same cofinality as $\theta$.
See the big picture of the construction on
page~\pageref{picture:construction}.  We intend to carry out  the
construction {\bf down} and {\bf to the right} in that picture.  In
the array, the bottom right hand corner ($M^*$) will be a
$(\mu,\omega)$-limit model witnessed by a chain of models as described in the
first paragraph of this section.  This chain will appear in the last
column of the array.  We will see that $M^*$ is a $(\mu,\theta)$-limit model
by examining the last (the $\omega^{\mbox{th}}$) row of the array. This last row will
be an $\prec_{\K}$-increasing sequence of models, $\bar M^*$ whose length
will have the same cofinality as $\theta$. However we will not be able to
guarantee that $M^*_{i+1}$ is universal over $M^*_i$ in this last row.
Thus we need another method to conclude that $M^*$ is a $(\mu,\theta)$-limit
model. This involves attaching more information to our sequence $\bar
M^*$.  We call this accessorized sequence of models a tower (see
Definition~\ref{def:towers} below).  Each row in our construction of
the array of models will be such a tower.

Under the assumption of Galois-superstability, given any sequence
$\langle a_i\mid
i<\theta\rangle$ of elements with $a_i\in M_{i+1}\backslash M_i$, we
can identify some
$N_i\prec_{\K}M_i$ such that $\tp(a_i/M_i)$ does not $\mu$-split over
$N_i$. Furthermore, by Assumption \ref{splitting assm}, we may choose
this $N_i$ such that $M_i$ is a limit model over $N_i$.  We abbreviate
this situation by the triple $(\bar M,\bar a,\bar N)$.

\begin{definition}[Towers]\label{def:towers}Let $(I,<)$ be a well ordering  of cardinality $<\mu^+$.
For cleaner notation, we will identify $I$ with $\theta$, its order-type, and we will denote
the successor of $i$ in the ordering $I$ by $i+1$ when it is clear.
Then, we define a \emph{tower} to be a triple $(\bar M,\bar a,\bar N)$
where $\bar M=\langle M_i\mid i<\theta\rangle$ is a $\prec_{\K}$-increasing sequence of limit
models of cardinality $\mu$; $\bar a=\langle a_i\mid i+1<\theta\rangle$ and $\bar N=\langle N_i\mid
i+1<\theta\rangle$ satisfy $a_i\in M_{i+1}\backslash M_i$; $\tp(a_i/M_i)$ does not $\mu$-split
over $N_i$; and $M_i$ is a $(\mu,\sigma)$-limit model over $N_i$. 
\end{definition}

\begin{notation}
We denote by $\K^*_{\mu,I}$ the set of towers of the form $(\bar
M,\bar a,\bar N)$ where 
the sequences $\bar M$, $\bar a$ and $\bar N$ are indexed by $I$.
Occasionally, $I$ will be an ordinal $\theta$ with the usual ordering,
and we write $\K^*_{\mu,\theta}$ for this set of towers.  
At times, we will be considering towers based on different well
orderings $I$ and $I'$ simultaneously.  In these contexts
 if $i\in I\bigcap I'$, the notation $i+1$ is not necessarily
 well-defined so we will use the
 notation $\s_I(i)$ for the successor of $i$ in the ordering $I$.
 Finally when $I$ is a sub-order of $I'$ for any $(\bar M,\bar a,\bar
N)\in\K^*_{\mu,I'}$ we write $(\bar M,\bar a,\bar N)\restriction I$ for
the tower in $\K^*_{\mu,I}$ given by the subsequences $\langle M_i\mid
i\in I\rangle$,
$\langle N_i\mid i\in I\rangle$ and $\langle a_i\mid i\in
I\rangle$.
\end{notation}

In addition to having control over the last row of the array, we also
need to be able to guarantee that  the last column of the tower
witnesses that $M^*$ is a $(\mu,\omega)$-limit model.  This will be done by
prescribing the following ordering on rows of the array:
\begin{definition}
For towers $(\bar M,\bar a,\bar N)\in\K^*_{\mu,I}$ and $(\bar M',\bar
a',\bar N')\in\K^*_{\mu,I'}$ with $I\subseteq I'$, we write
$(\bar M,\bar a,\bar N)<(\bar M',\bar a',\bar N')$ if and only if
for every $i\in I$, $a_i=a'_i$, $N_i=N'_i$ and $M'_i$ is a proper
universal extension of $M_i$.

\end{definition}

\begin{remark} The ordering $<$ on towers is identical to the ordering
$<^c_\mu$ defined in \cite{ShVi}.  The superscript was used by Shelah and
Villaveces to distinguish this ordering from others.  We only use one
ordering on towers, so we omit the superscripts and subscripts here.
\end{remark}

Once we have established an ordering on towers, we can define a
specific tower which will be called a {\em union of an increasing
  sequence of towers}. Suppose that $\langle(\bar M,\bar a,\bar
N)^\gamma\in\K^*_{\mu,I_\gamma}\mid \gamma<\beta\rangle$ is an
increasing sequence of towers such that
the index set $I_\gamma$ of $(\bar M,\bar a,\bar N)^\gamma$ is a
sub-ordering of
the index set $I_{\gamma'}$ for $(\bar M,\bar a,\bar N)^{\gamma'}$
whenever
$\gamma<\gamma'$. Let  $I_\beta:=\Union_{\gamma<\beta}I_\gamma.$  Then
denote by $(\bar M,\bar
a,\bar N)^\beta\in\K^*_{\mu,I_\beta}$ the ``union'' of the sequence of
towers where
\[ a^\beta_i=a^{\min\{\gamma\mid i\in I_\gamma\}}_i,\]
\[ N^\beta_i= N^{\min\{\gamma\mid i\in I_\gamma\}}_i
\text{ and}\]
\[\bar M^\beta=\langle M^\beta_i\mid i\in \Union_{\gamma<\beta}I_\gamma\rangle
\text{ with }
M^\beta_i=\bigcup_{\gamma<\beta}\bigcup_{I_\gamma \ni i}M^\gamma_i.\]
By Assumption \ref{splitting assm}.\ref{locality}, $(\bar M,\bar
a,\bar N)^\beta$ is indeed a tower.

Notice that we do not assume an individual tower to be continuous.
Nor do we assume that inside of a tower $M_{i+1}$ is universal over
$M_i$.  If one considers the approach of defining an array of models
row by row, then generally (even in the first order case)  even if all
rows are continuous and satisfy the universality property mentioned in
this paragraph, it is not necessarily true that the union of these
rows will be a tower in which every model is universal over its
predecessors.

%

For a tower $(\bar M,\bar a,\bar N)$, it was shown in \cite{ShVi},
that even if $M_{i+1}$ is not universal over $M_i$, one can conclude
that $\Union_{i<\theta}M_i$ is a $(\mu,\theta)$-limit model provided that all
types over each of  the $M_i$ are realized by a sufficient number of
$a_j$s in the tower.  Unfortunately constructing such a tower meeting
these along with all of our other requirements is beyond
reach. However, in \cite{Va1}, VanDieren showed that slightly less was
needed (see Definition \ref{def:relativefulltowers}).  
In \cite{Va1}, the amalgamation property is not assumed resulting in
noise that can be avoided in our context.  Thus because we have at our
disposal the AP, we provide 
a complete, undistracted proof here.

\begin{definition}[Relatively Full Towers]\label{def:relativefulltowers}
  Suppose that $I$ is a well-ordered set such that there exists a
  cofinal sequence $\langle i_\alpha\mid\alpha<\theta\rangle$ of $I$
  of order type $\theta$ such that
  there are $\mu\cdot \omega$ many elements between $i_\alpha$ and
  $i_{\alpha+1}$.\\
  Let $(\bar M,\bar a,\bar N)$ be a tower indexed by $I$ such that
  each $M_i$ is a $(\mu,\sigma)$-limit model.  For each
$i$, let
$\langle M^\gamma_{i}\mid \gamma<\sigma\rangle$ witness that
$M_{i}$ is a
$(\mu,\sigma)$-limit model.\\
The tower
$(\bar M,\bar a,\bar N)$ is
\emph{full relative to
$(M^\gamma_{i})_{\gamma<\sigma,i\in I}$} iff for every
$\gamma<\sigma$ and every $(p,M^\gamma_{i})\in\St(M_{i})$ with
$i_\alpha\leq i<i_{\alpha+1}$, there exists
$j\in I$   with $i\leq j< i_{\alpha+1}$ such that
$(\tp(a_j/M_j),N_j)$ and
$(p,M^\gamma_{i})$ are parallel.
\end{definition}

Relative fullness of towers can be seen as a (weak) form of ``eventual
Galois saturation.'' Along a full tower, all strong Galois types over
members of sequences - sequences which witness the fact that along the tower the
models are limits - end up being realized (modulo parallelism) by an
element $a_j$ of the tower. As we see in our proof, this property is
much more flexible than regular Galois-saturation - it could be
regarded as a ``dynamic'' and robust version.

Although relatively full towers are used here as a technical device
for the proof, the crucial property is that ``eventual'' or
``dynamic'' relative Galois saturation. These objects have variously
been used by Shelah in various places, Shelah-Villaveces~\cite{ShVi},
VanDieren~\cite{Va1}, and other authors. It is reasonable to say that
the notion of relatively full towers has potential for other uses
outside of these works.

\begin{theorem}[Relatively full towers provide limit
  models]\label{relatively full is limit} Let $\theta$ be a limit ordinal
  $<\mu^+$ satisfying $\theta=\mu\cdot\theta$.  Suppose that $I$ is a
  well-ordered set
  as in Definition \ref{def:relativefulltowers}.

Let $(\bar M,\bar a,\bar N)\in\K_{\mu,I}^*$ be a tower made up of
$(\mu,\sigma)$-limit models, for some fixed $\sigma<\mu^+$. If $(\bar M,\bar a,\bar
N)\in\K^*_{\mu,I}$ is full relative to $(M^\gamma_i)_{i\in  I,\gamma<\sigma}$, then
$M:=\Union_{i\in I}M_i$ is a $(\mu,\theta)$-limit model.
\end{theorem}

\begin{proof}
  Without loss of generality we may assume that $\bar M$ is
  continuous.
Let $M'$ be a $(\mu,\theta)$-limit model over $M_{i_0}$ witnessed by
$\langle M'_\alpha\mid\alpha<\theta\rangle$.  By $\mu$-DAP over limit models, we may assume that $M'\cap M=M_{i_0}$.  Since
$\theta=\mu\cdot\theta$, we may also arrange things so that the universe
of $M'_\alpha$ is $\mu\cdot \alpha$ and $\alpha\in M'_{\alpha+1}$.

 We will construct an
isomorphism between $M$ and $M'$ by induction on $\alpha<\theta$.
Define an increasing and
continuous sequence of 
$\prec_{\K}$-mappings
$\langle h_\alpha\mid \alpha<\theta\rangle$ such that
\begin{enumerate}
\item\label{hi cond} $h_\alpha:M_{i_\alpha+j}\to M'_{\alpha+1}$
for some
$j<\mu\cdot\omega$
\item $h_0=\id_{M_{i_0}}$ and
\item\label{put it all in} $\alpha\in \rg(h_{\alpha+1})$.
\end{enumerate}
For $\alpha=0$ take $h_0=\id_{M_{i_0}}$.  For $\alpha$ a limit ordinal let
$ h_\alpha=\Union_{\beta<\alpha}h_{\beta}$.  Since $\bar M$ is
continuous, the induction hypothesis gives us that $h_\alpha$ is a
$\prec_{\K}$-mapping from
$M_{i_\alpha}$ into
$M'_\alpha$ allowing us to satisfy condition (\ref{hi cond}) of the
construction.

Suppose that $h_\alpha$
has been defined.  Let $j<\mu\cdot\omega$ be such that
$h_\alpha:M_{i_\alpha+j}\to M'_{\alpha+1}$.   There are two
cases:  either $\alpha\in\rg(h_\alpha)$ or
$\alpha\notin\rg(h_\alpha)$. First suppose that $\alpha\in\rg(h_\alpha)$. 
Since
$M'_{\alpha+2}$ is universal over
$M'_{\alpha+1}$, it is also universal over $h_\alpha(M_{i_\alpha+j})$. 
This allows us to extend $h_\alpha$ to
$h_{\alpha+1}:M_{i_{\alpha+1}}\to M'_{\alpha+2}$.

Now consider the case when $\alpha\notin\rg(h_\alpha)$.
Since $\langle M^\gamma_{i_\alpha+j}\mid \gamma<\sigma\rangle$ witnesses
that $M_{i_\alpha+j}$ is a $(\mu,\sigma)$-limit model, by Assumption
\ref{splitting assm}, there exists
$\gamma <\sigma$ such that
$\tp(\alpha/M_{i_\alpha+j})$ does not $\mu$-split over
$M^\gamma_{i_\alpha+j}$.
By our choice of $\bar M'$ disjoint from $\bar M$ outside of
$M_{i_0}$, we know that $\alpha\notin M_{i_\alpha+j}$.  Thus
$\tp(\alpha/M_{i_\alpha+j})$ is non-algebraic.  By relative fullness of
$(\bar M,\bar a,\bar N)$, there exists
$j'$ with $j\leq j'<i_{\alpha+1}$ such that
$(\tp(\alpha/M_{i_\alpha+j'}),M^\gamma_{i_\alpha+j})$
is parallel to $(\tp(a_{i_{\alpha+1}+j'}/M_{i_{\alpha+1}+j'}),
N_{i_{\alpha+1}+j'})$.  In particular we have
that
$$(*)\quad\tp(a_{i_{\alpha+1}+j'}/M_{i_\alpha+j})=\tp(\alpha/M_{i_\alpha+j}).$$

We can extend $h_\alpha$ to an automorphism $h'$ of $\C$.
An application of $h'$ to $(*)$ gives us

$$(**)\quad\tp(h'(a_{i_{\alpha+1}+j'})/h_\alpha(M_{i_\alpha+j}))=
\tp(\alpha/h_\alpha(M_{i_\alpha+j})).$$

Since $M'_{\alpha+2}$ is universal over $h_{\alpha}(M_{i_\alpha})$, we may
extend $h_{\alpha}$ to a $\K$-mapping
$h_{\alpha+1}:M_{i_{\alpha+1}+j'}\to M'_{\alpha+2}$ such that
$h_{\alpha+1}(a_{i_{\alpha+1}+j'})=\alpha$.

Let $h:=\Union_{\alpha<\theta}h_\alpha$.  Clearly
$h:M\to M'$.  To see that $h$ is an
isomorphism, notice that condition (\ref{put it all in}) of the
construction forces $h$ to be surjective.
\end{proof}

\section{Uniqueness of Limit Models}

We now begin the construction of our array of models and $M^*$.
Let $\theta$ be an ordinal as in the previous section.  The
goal is to build an array of models with $\omega+1$ rows so that the
bottom row of the array is a relatively full tower indexed by a set of cofinality $\theta$. 
 To do this, we will be
adding elements to the index set of towers row by row so that 
at stage $n$ of our construction the
tower that we build is indexed by $I_n$
 described here:

\begin{notation}
The index sets $I_n$ will be defined inductively so that $\langle I_n\mid
n<\omega+1\rangle$ is an increasing and continuous chain of well-ordered sets. We
fix $I_0$ to be an index set of order type $\theta+1$ and will denote it by
$\langle i_\alpha\mid\alpha\leq\theta\rangle$.  We will refer to the members of $I_0$ by name in many
stages of the construction.  These indices serve as anchors for the
members of the remaining index sets in the array.  Next we demand that
for each $n<\omega$, $\{j\in I_n\mid i_\alpha<j<i_{\alpha+1}\}$ has order type $\mu\cdot n$ such
that each $I_n$ has supremum $i_\theta$. An example of such $\langle I_n\mid n\leq\omega\rangle$
is $I_n=\theta\times(\mu\cdot n)\Union\{i_\theta\}$ ordered lexicographically, where $i_\theta$ is
an element $\geq$  each $i\in \Union_{n<\omega}I_n$. Also, let $I=\bigcup_{n<\omega}I_n$.
\end{notation}

To prove the main theorem of the paper, we need to prove that for a
fixed $M\in\K$ of cardinality $\mu$ any $(\mu,\theta)$-limit and $(\mu,\omega)$-limit model
over $M$ are isomorphic over $M$.  Let  us begin by fixing $M\in\K_\mu$ and $\theta$ such that
$\mu\cdot\theta=\theta$.  Without loss of generality, $M$ is a limit model.
We define by induction on $n\leq\omega$ a $<$-increasing and continuous
sequence of towers $(\bar M,\bar a,\bar N)$ such that
\begin{enumerate}
\item $(\bar M,\bar a,\bar N)^0$ is a tower with $M^0_0=M$.
\item $(\bar M,\bar a,\bar N)\in\K^*_{\mu,I_n}$
\item For every $(p,N)\in\St(M^n_i)$ with $i_\alpha\leq i< i_{\alpha+1}$ there is 
$j\in I_{n+1}$ with $i<j<i_{\alpha+1}$ so that $(\tp(a_j/M^{n+1}_j),N^{n+1}_j)$ and $(p,N)$ 
are parallel.
\item $M^{n+1}_{i_{\alpha+1}}$ is a $(\mu,\mu)$-limit
model over $\bigcup_{j<i_{\alpha+1}}M^{n+1}_j$.

\end{enumerate}

Given $M$, we can find a tower $(\bar M,\bar a,\bar N)^0\in\K^*_{\mu,I_0}$
with $M^0_0=M$ because of the existence of universal extensions and
because of Assumption \ref{splitting assm}.\ref{no long splitting
  chain}. The last pages (Page \pageref{picture:construction} onward)
of this section provide a picture of this construction of an array of
models, explanations for carrying out the final stage of the
construction and a proof that this is sufficient to prove the main
theorem. We spend most of the remainder of this section verifying
that it is possible to carry out the induction step of the
construction.  This is a particular case of Theorem II.7.1 of
\cite{Va1}. But since our context is somewhat easier, we do not
encounter so many obstacles as in \cite{Va1} and we provide a
different, more direct proof here:

\begin{theorem}[Dense $<$-extension property]\label{<c}
  Given $(\bar M,\bar a,\bar N)\in \K^*_{\mu, I_{n}}$ there exists $(\bar
  M',\bar a,\bar N)\in\K^*_{\mu, I_{n+1}}$ such that $(\bar M,\bar a,\bar
  N)<(\bar M',\bar a,\bar N)$ and for each $(p,N)\in\St(M_{i})$ with
  $i_\alpha\leq i<i_{\alpha+1}$, there exists $j\in I_{n+1}$ with $i<j<i_{\alpha+1}$ such
  that $(\tp(a_j/M'_j),N_j)$ and $(p,N)$ are parallel. Here,  the
  $M_i$'s are defined for $i\in I_n$ and the $M'_j$ are defined for $j\in
  I_{n+1}$.
\end{theorem}

Before we prove Theorem \ref{<c}, we prove a slightly weaker extension
property, one in which we can find an extension of the tower $(\bar M,
\bar a, \bar N)$ of the same index set:

\begin{lemma}[$<$-extension property]\label{<c-lemma}
  Given $(\bar M,\bar a,\bar N)\in\K^*_{\mu,I}$,
  there exists a $<$-extension $(\bar M',\bar a,\bar N)\in\K^*_{\mu,I}$ of
  $(\bar M,\bar a,\bar N)$ such that for each
limit $i$, $M'_i$ is a $(\mu,\mu)$-limit model over $\bigcup_{j<i}M'_j$.
\end{lemma}
\begin{proof}


  Given $(\bar M,\bar a,\bar N)\in\K^*_{\mu,I}$ we will define a
  $<$-extension $(\bar M',\bar a,\bar N)$ by induction on $i\in
  I$. Notice that a straightforward induction proof is not sufficient
  here for if we have defined $\langle M_j\mid j\leq i\rangle$ as a tower extending
  $(\bar M,\bar a,\bar N)$ restricted to $\langle j\mid j\leq i\rangle$ and are at the
  stage of defining $M'_{i+1}$, we may be faced with an impossible
  task: during our construction we may have inadvertently placed
  inside $M'_i$ witnesses for the splitting of the type of $a_{i+1}$
  over $N_{i+1}$; this would prevent us from extending $M'_i$ to
  $M'_{i+1}$ so that $\tp(a_{i+1}/M'_{i+1})$ does not $\mu$-split over
  $N_{i+1}$.  Therefore, we will instead define approximations, $M^+_i$, for $M'_i$ by induction on $i\in I$ and at each stage $i$ of
  the induction we will make adjustments of the previously defined
  approximation $M^+_j$ for $j<i$.  This leads us into defining
  $M^+_i$ and a directed system of $\prec_{\K}$-embeddings $\langle f_{j,i}\mid
  j<i\in I\rangle$ such that for $i\in I$, $M_i \prec_{\K} M^+_i$ for $j\leq i$,
  $f_{j,i}:M^+_j\to M^+_i$ and $f_{j,i}\restriction M_j=\id_{M_j}$.  We
  further require that $M^+_{i+1}$ is a limit model over
  $f_{i,i+1}(M^+_i)$ and $\tp(a_i/f_{i,i+1}(M^+_i))$ does not
  $\mu$-split over $N_i$.  When $i$ is a limit, we choose $M^+_i$ to be
  a $(\mu,\mu)$-limit model over $\Union_{j<i}f_{j,i}(M^+_j)$.

  This construction is done by induction on $i\in I$ using the existence
  of non-$\mu$-splitting extensions.
Suppose that $\langle M^+_k\mid k\leq i\rangle$ and $\langle f_{k,l}\mid k\leq l\leq i\rangle$ have been defined.
We explain how to define $M^+_{i+1}$ and $f_{i,i+1}$.  The rest of the
definitions required for the $i+1^{\mbox{st}}$ stage are dictated by
the requirement that we are forming a directed system. Let $M^*_{i+1}$
be an  limit model over both $M^+_i$ and  $M_{i+1}$.
Since $\tp(a_{i+1}/M_{i+1})$ does not $\mu$-split over $N_{i+1}$, 
by Fact \ref{splitting extension lemma} there exists
$f\in\Aut_{M_{i+1}}(\C)$ so that $\tp(a_{i+1}/f(M^*_{i+1}))$ does not
$\mu$-split over $N_{i+1}$.  Take $M^+_{i+1}:=f(M^*_{i+1})$ and
$f_{i,i+1}:=f\restriction M^+_i$.


At limit stages we take direct limits so that $f_{j,i}\restriction
M_j=\id_{M_j}$.  This is possible by Subclaims II.7.10 and II.7.11 of
\cite{Va1} or see Claim 2.17 of
\cite{GrVa2}.  Take an extension of the direct limit that is both
universal over $M_i$ and
is a $(\mu,\mu)$-limit over $\bigcup_{j< i} f_{j,i}(M_j)$ and call this
$M^+_i$. Notice that we do not obtain a continuous tower; 
continuity will be recovered later using reduced towers.

Let $f_{j,\sup\{I\}}$ and $M'_{\sup\{I\}}$ be the direct limit of this
system such that $f_{j,\sup\{I\}}\restriction M_j=\id_{M_j}$. We can now
define
$M'_j:=f_{j,\sup\{I\}}(M^+_j)$ for each $j\in I$.
By construction, we have that $\tp(a_i/f_{i,i+1}(M_i^+))$ does not
$\mu$-split over $N_i$. Mapping into $M_{\sup(I)}$ by $f_{i+1,\sup(I)}$,
and noting that both $a_i$ and $N_i$ are fixed by $f_{i+1,\sup(I)}$,
we conclude that $\tp(a_i/M_i')$ does not $\mu$-split over $N_i$ as required.


\end{proof}

We can now use the extension property for towers of the same index set
from Lemma \ref{<c-lemma} to prove the dense extension property which
allows us to grow the index set as we add elements to the models in
the extension.

\begin{proof}[Proof of Theorem \ref{<c}]

Given $(\bar M,\bar a,\bar N)\in\K^*_{\mu,I_n}$, let $(\bar M',\bar
a,\bar N)\in\K^*_{\mu,I_n}$ be an extension of $(\bar M,\bar a,\bar N)$
as in Lemma \ref{<c-lemma} so that each $M'_{i_{\alpha+1}}$ is a
$(\mu,\mu)$-limit model over $\bigcup_{j<i_{\alpha+1}}M'_j$.

For each $i_\alpha$, let $\langle M'_l\mid l\in I_{n+1},\;i_\alpha+\mu\cdot n<l<i_{\alpha+1}\rangle$ witness
that $M'_{i_{\alpha+1}}$ is a $(\mu,\mu)$-limit model over
$\bigcup_{j<i_{\alpha+1}}M'_{j}$.  Without loss of generality we may assume that
each of these $M'_l$ is a limit model over its predecessor.  

Fix $\{(p,N)^l_{i_\alpha}\mid i_\alpha+\mu\cdot n<l<i_{\alpha+1}\}$ an enumeration of
  $\Union\{ \St(M_i) : i\in I_n,i_\alpha\leq i<i_{\alpha+1} \}$. By our choice of
  $I_{n+1}$ and stability in $\mu$, such an enumeration is possible.
  Since  $M'_{\s_{I_{n+1}}(l)}$ is universal over $M'_l$, there exists
  a realization in $M'_{\s_{I_{n+1}}(l)}$ of the non-$\mu$-splitting
  extension of $p_{i_\alpha}^l$ to $M'_l$.  Let $a_l$ be this realization
  and take $N_l:=N_{i_\alpha}^l$.

Notice that $(\langle M'_j\mid j\in I_{n+1}\rangle,\langle a_j\mid j\in I_{n+1}\rangle,\langle N_j\mid j\in
I_{n+1}\rangle)$ provide the desired extension of $(\bar M,\bar a,\bar N)$
in $\K^*_{\mu,I_{n+1}}$.
\end{proof}

We are almost ready to carry out the complete construction.  However,
notice that Theorem \ref{<c} does not provide us with a continuous
extension.  Therefore the bottom (i.e. the $\omega+1^{st}$) row of our
array may not be continuous
which would prevent us from applying Theorem~\ref{relatively
full is limit} to conclude that $M^*$ is a $(\mu,\theta)$-limit model. 
So we will further require that the towers that occur in the rows of
our array are all
continuous.  This can be guaranteed by restricting ourselves to reduced
towers as in \cite{ShVi} and \cite{Va1}.
\begin{definition}\label{reduced defn}\index{reduced towers}
A tower $(\bar M,\bar a,\bar N)\in\K^*_{\mu,I}$ is said to 
be \emph{reduced} provided that for every $(\bar M',\bar a,\bar
N)\in\K^*_{\mu,I}$ with
$(\bar M,\bar a,\bar N)\leq(\bar M',\bar a,\bar
N)$ we have that for every
$i\in I$,
$$(*)_i\quad M'_i\cap\Union_{j\in I}M_j = M_i.$$
\end{definition}

If we take a $<$-increasing
chain of reduced towers, the union
will be reduced.  The following fact appears as
Theorem 3.1.14 of \cite{ShVi}.  We provide the
proof for completeness.

\begin{fact}\label{union of reduced is reduced}
Let
$\langle (\bar M,\bar a,\bar N)^\gamma\in\K^*_{\mu,I_\gamma}\mid
\gamma<\beta\rangle$ be a $<$-increasing and continuous
sequence of
reduced towers such that the sequence is continuous in the sense that
for a limit $\gamma<\beta$, the tower $(\bar M,\bar a,\bar N)^\gamma$ is the union of the
towers $(\bar M,\bar a,\bar N)^\zeta$ for $\zeta<\gamma$.
Then the union of the sequence of towers $\langle (\bar M,\bar a,\bar N)^\gamma\in\K^*_{\mu,I_\gamma}\mid
\gamma<\beta\rangle$ is itself a
reduced tower.
\end{fact}
\begin{proof}

  Suppose that $(\bar M,\bar a,\bar N)^\beta$ is not reduced. 
Let
$(\bar M',\bar a,\bar N)\in\K^*_{\mu,I_\beta}$ witness this.
Then there exists an $i\in I_\beta$ and an element $b$ such that
$b\in (M'_i\cap\Union_{j\in I_\beta}M^\beta_j)\backslash M^\beta_i$.
There exists $\gamma<\beta$ such that 
$b\in \Union_{j\in I_\gamma}M^\gamma_j\backslash M^\gamma_i$.  
Notice that $(\bar M',\bar a,\bar N)\restriction I_\gamma$ witnesses that
$(\bar M,\bar a,\bar N)^\gamma$ is not reduced.
\end{proof}

The following 
appears in \cite{ShVi} (Theorem 3.1.13).

\begin{fact}[Density of reduced towers]\label{density of
reduced}\index{reduced towers!density of} There exists a reduced
$<$-extension of every tower in
$\K^*_{\mu,I}$.

\end{fact}
\begin{proof}

Assume for the sake of contradiction that no $<$-extension of $(\bar
M,\bar a,\bar N)$ is reduced.  This allows us to construct a
$\leq$-increasing and continuous sequence of towers $\langle(\bar M,\bar a,\bar
N)^\zeta\in\K^*_{\mu,I}\mid\zeta<\mu^+\rangle$ such that $(\bar M,\bar a,\bar N)^{\zeta+1}$
witnesses that $(\bar M,\bar a,\bar N)^\zeta$ is not reduced. The
construction is done inductively in the obvious way.

%
%

For each $b\in\Union_{\zeta<\mu^+,i\in I}M^\zeta_i$ define 
$$i(b):=\min\big\{i\in I\mid
b\in\Union_{\zeta<\mu^+}\Union_{j\leq i}
M^\zeta_j\big\}\text{ and}$$
$$\zeta(b):=\min\big\{\zeta<\mu^+\mid
b\in M^\zeta_{i(b)}\big\}.$$

$\zeta(\cdot)$ can be viewed as a function from $\mu^+$ to $\mu^+$.  
Since $|I|=\mu$ and each $M^\zeta_i$ has cardinality $\mu$,
there exists a club $E=\{\delta<\mu^+\mid \forall b\in
\Union_{i\in I}M^\delta_i,\; \zeta(b)<\delta\}$.  Actually, all we
need is for $E$ to be non-empty.

Fix $\delta\in E$.  By construction
$(\bar M,\bar
a,\bar N)^{\delta+1}$
witnesses the fact that $(\bar M,\bar a,\bar
N)^\delta$ is not reduced.  
So we may fix $i\in I$ and $b\in
M^{\delta+1}_i\cap\Union_{j\in I}M^{\delta}_j$ such that $b\notin
M^\delta_i$.  Since $b\in M^{\delta+1}_i$, we have that $i(b)\leq i$.
Since $\delta\in E$, we know that there exists $\zeta<\delta$ such that
$b\in M^{\zeta}_{i(b)}$.  Because $\zeta<\delta$ and $i(b)\leq i$, this
implies that $b\in M^\delta_i$ as well.  This contradicts our choice of
$i$ and $b$ witnessing the failure of $(\bar M,\bar a,\bar
N)^\delta$ to be reduced.
\end{proof}

By revising the proof of Lemma \ref{<c-lemma}, we can conclude:
\begin{lemma}\label{monotonicity of C-red}
Suppose that $(\bar M,\bar a,\bar N)\in\K^*_{\mu,I}$ is
reduced.  If  $I_0$ is an initial segment of $I$, then $(\bar M,\bar
a,\bar N)\restriction I_0$ is reduced.
\end{lemma}

\begin{proof}
Suppose that $(\bar M,\bar a,\bar N)\restriction I_0$ is not reduced.  
Let $(\bar M',\bar a\restriction I_0,\bar N\restriction I_0)$ and
$\delta<j\in I_0$ with $b\in (M_\delta'\cap M_j)\backslash M_\delta$ witness this.
We can apply the inductive step of
Lemma \ref{<c-lemma} (replacing an initial segment of the construction
there with $\bar{M}'$), to find $(\bar M'',\bar a,\bar N)$ an extension
of $(\bar M,\bar a,\bar N)$ such that there is a $\prec_{\K}$-mapping $f$
from the models of $\bar M'$ into the models of $\bar M''$ with
$f\restriction M_j=\id_{M_j}$.  Notice that $(\bar M'',\bar a,\bar N)$
and $b,\delta,j$ will witness that $(\bar M,\bar a,\bar N)$ is not reduced.
\end{proof}

The following theorem makes use of the unidimensionality assumption.
This generalizes a special case of the uniqueness of limit models
result in the series of papers \cite{Va1} and \cite{Va2} by replacing
the assumption of categoricity in $\mu^+$ with the weaker
unidimensionality assumption.  Further work of VanDieren in \cite{Va3}
weakens this assumption further for tame classes.

\begin{theorem}[Reduced towers are
continuous]\label{reduced are cont}\index{reduced towers!are
continuous}If $(\bar M,\bar a,\bar N)\in\K^*_{\mu,I}$ is reduced, then it
is continuous, namely for each limit $i$ in $I$, $M_i=\bigcup_{j<i}M_j$.

\end{theorem}

\begin{proof}[Proof of Theorem \ref{reduced are cont}]

Suppose the theorem fails for $\mu$.
Let  $\delta$ be the minimal limit ordinal such that there exists an index set $I$
and  $(\bar M,\bar a,\bar N)\in\K^*_{\mu,I}$ a reduced tower which is discontinuous
at the $\delta^{\mbox{th}}$ element of $I$.  We can apply Lemma
\ref{monotonicity of C-red} to assume without loss of generality that
$I=\delta+1$.
Fix $(\bar M,\bar a,\bar N)\in\K^*_{\mu,\delta+1}$  reduced and discontinuous
at $\delta$ with $b\in M_\delta\backslash\Union_{i<\delta}M_i$.  By
Fact \ref{existence of minimal}, there exists a minimal type $p$ over
$M_0$.  So by our unidimensionality Assumption
\ref{unidimensionality}, we 
know that the Galois type of $p$ must be realized in $M_\delta\backslash\Union_{i<\delta}M_i$.  Therefore, we may  
assume that $b\models p$.

\begin{claim}\label{b in claim}
There exists a $<$-extension of $(\bar M,\bar a,\bar
N)\restriction\delta$, containing $b$.  We will refer to such a tower
in $\K^*_{\mu,\delta}$ as
$(\bar M',\bar a\restriction\delta,\bar
N\restriction\delta)$. Furthermore, $b$ may be assumed to be an
element of $M_0'$.

\end{claim}

\begin{proof}[Proof of Claim \ref{b in claim}]
We use the minimality of $\delta$ and the $<$-extension property to find a tower of length $\delta$,  $(\bar M^*,\bar a\restriction\delta,\bar N\restriction\delta)$, that 
is a proper extension of $(\bar M,\bar a,\bar N)\restriction\delta$.  By the definition of $<$-extension, $M^*_0$ is universal over $M_0$; so we can find $b^*\in M^*_0\backslash M_0$ realizing $p$.  

Notice that by Lemma \ref{monotonicity of C-red}, $(\bar M,\bar a,\bar N)\restriction\delta$ is reduced.  Thus we can conclude that $b^*\in M^*_0\backslash \Union_{i<\delta}M_i$
and $\tp(b^*/\Union_{i<\delta}M_i)$ is non-algebraic.  Since $p$ is minimal, it must be the case that $\tp(b^*/\Union_{i<\delta}M_i)=\tp(b/\Union_{i<\delta}M_i)$.  Let $f\in\Aut_{\Union_{i<\delta}M_i}\C$ take $b^*$ to $b$.  

Consider the image of $(\bar M^*,\bar a,\bar N)$ under $f$; denote this tower by $(\bar M',\bar a,\bar N)$.  Because $f$ fixes $(\bar M,\bar a,\bar N)\restriction \delta$, $(\bar M',\bar a,\bar N)$ is an extension of $(\bar M,\bar a,\bar N)\restriction \delta$ as required.  
\end{proof}

Using $(\bar M',\bar a,\bar N)$ from Claim \ref{b in claim}, define $M'_\delta$ to be a limit model of cardinality $\mu$ containing $\Union_{i<\delta}M'_i$ so that it is universal over $M_\delta$.  Notice that the tower $(\bar M'\conc \langle M'_\delta\rangle,\bar a,\bar N)$ extends $(\bar M,\bar a,\bar N)$ with $b\in (M'_0\backslash \Union_{i<\delta}M_i)\bigcap M_\delta$.  This contradicts our assumption that $(\bar M,\bar a,\bar N)$ is reduced and completes the proof of Theorem \ref{reduced are cont}.

\end{proof}

\begin{corollary}\label{continuous extension}
In Theorem \ref{<c}, we can choose $(\bar M,\bar a,\bar N)$ to 
be reduced, and hence continuous.
\end{corollary}

Now we return to the construction in the proof of the Main Theorem.

\medskip

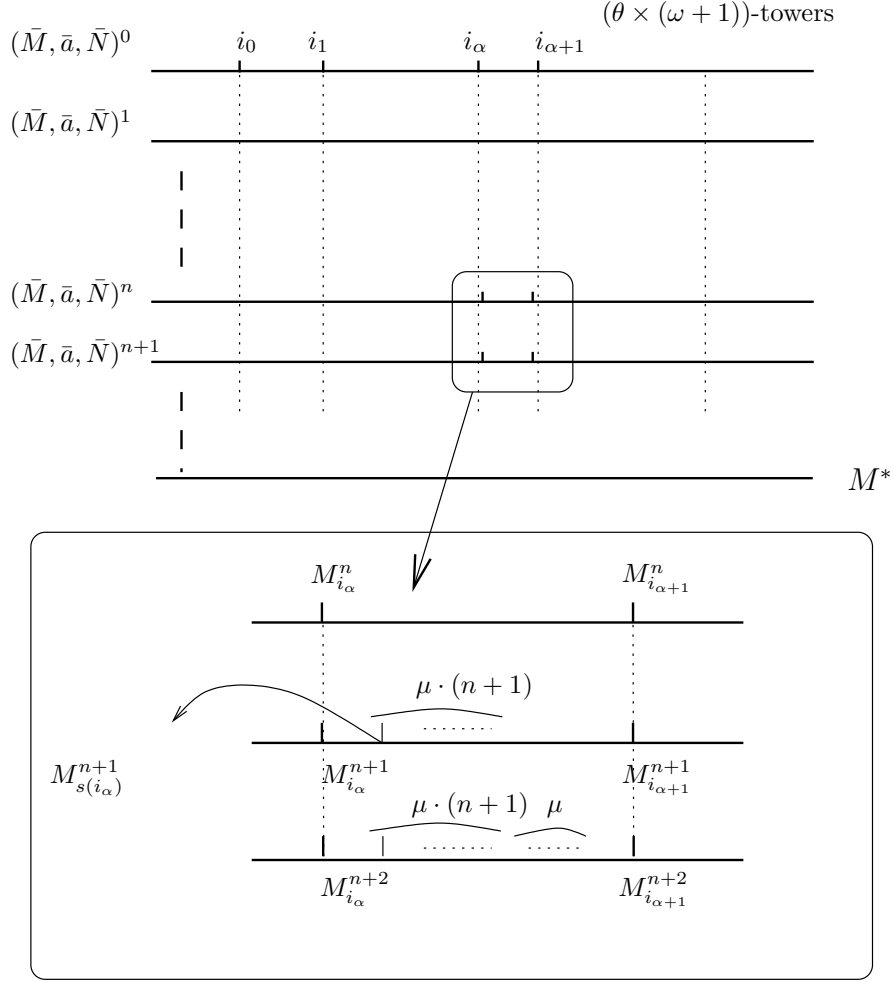
\begin{picture}(0,0)%
\includegraphics{GVVb.pstex}%
\end{picture}%
\setlength{\unitlength}{3947sp}%
\begingroup\makeatletter\ifx\SetFigFont\undefined%
\gdef\SetFigFont#1#2#3#4#5{%
  \reset@font\fontsize{#1}{#2pt}%
  \fontfamily{#3}\fontseries{#4}\fontshape{#5}%
  \selectfont}%
\fi\endgroup%
\begin{picture}(5711,6156)(76,-5623)
\put( 76,-1352){\makebox(0,0)[lb]{\smash{{\SetFigFont{10}{12.0}{\familydefault}{\mddefault}{\updefault}{ $(\bar{M},\bar{a},\bar{N})^n$}%
}}}}
\put( 76,-1730){\makebox(0,0)[lb]{\smash{{\SetFigFont{10}{12.0}{\familydefault}{\mddefault}{\updefault}{ $(\bar{M},\bar{a},\bar{N})^{n+1}$}%
}}}}
\put(3921,-3116){\makebox(0,0)[lb]{\smash{{\SetFigFont{10}{12.0}{\familydefault}{\mddefault}{\updefault}{ $M^n_{i_{\alpha+1}}$}%
}}}}
\put(2030,-4377){\makebox(0,0)[lb]{\smash{{\SetFigFont{10}{12.0}{\familydefault}{\mddefault}{\updefault}{ $M_{i_\alpha}^{n+1}$}%
}}}}
\put(3921,-4377){\makebox(0,0)[lb]{\smash{{\SetFigFont{10}{12.0}{\familydefault}{\mddefault}{\updefault}{ $M_{i_{\alpha +1}}^{n+1}$}%
}}}}
\put(328,-4377){\makebox(0,0)[lb]{\smash{{\SetFigFont{10}{12.0}{\familydefault}{\mddefault}{\updefault}{ $M^{n+1}_{s(i_\alpha)}$}%
}}}}
\put(3451,-4561){\makebox(0,0)[lb]{\smash{{\SetFigFont{10}{12.0}{\familydefault}{\mddefault}{\updefault}{ $\mu$}%
}}}}
\put(2026,-5086){\makebox(0,0)[lb]{\smash{{\SetFigFont{10}{12.0}{\familydefault}{\mddefault}{\updefault}{ $M_{i_\alpha}^{n+2}$}%
}}}}
\put(3901,-5086){\makebox(0,0)[lb]{\smash{{\SetFigFont{10}{12.0}{\familydefault}{\mddefault}{\updefault}{ $M_{i_{\alpha +1}}^{n+2}$}%
}}}}
\put( 76,224){\makebox(0,0)[lb]{\smash{{\SetFigFont{10}{12.0}{\familydefault}{\mddefault}{\updefault}{ $(\bar{M},\bar{a},\bar{N})^0$}%
}}}}
\put(1501,239){\makebox(0,0)[lb]{\smash{{\SetFigFont{10}{12.0}{\familydefault}{\mddefault}{\updefault}{ $i_0$}%
}}}}
\put(1951,239){\makebox(0,0)[lb]{\smash{{\SetFigFont{10}{12.0}{\familydefault}{\mddefault}{\updefault}{ $i_1$}%
}}}}
\put( 76,-280){\makebox(0,0)[lb]{\smash{{\SetFigFont{10}{12.0}{\familydefault}{\mddefault}{\updefault}{ $(\bar{M},\bar{a},\bar{N})^1$}%
}}}}
\put(2926,239){\makebox(0,0)[lb]{\smash{{\SetFigFont{10}{12.0}{\familydefault}{\mddefault}{\updefault}{ $i_\alpha$}%
}}}}
\put(3376,239){\makebox(0,0)[lb]{\smash{{\SetFigFont{10}{12.0}{\familydefault}{\mddefault}{\updefault}{ $i_{\alpha+1}$}%
}}}}
\put(5326,-2536){\makebox(0,0)[lb]{\smash{{\SetFigFont{12}{14.4}{\familydefault}{\mddefault}{\updefault}{ $M^*$}%
}}}}
\put(3795,413){\makebox(0,0)[lb]{\smash{{\SetFigFont{10}{12.0}{\familydefault}{\mddefault}{\updefault}{ $(\theta \times (\omega+1))$-towers}%
}}}}
\put(1967,-3116){\makebox(0,0)[lb]{\smash{{\SetFigFont{10}{12.0}{\familydefault}{\mddefault}{\updefault}{ $M_{i_\alpha}^n$}%
}}}}
\put(2626,-3811){\makebox(0,0)[lb]{\smash{{\SetFigFont{10}{12.0}{\familydefault}{\mddefault}{\updefault}{ $\mu \cdot (n+1)$}%
}}}}
\put(2600,-4561){\makebox(0,0)[lb]{\smash{{\SetFigFont{10}{12.0}{\familydefault}{\mddefault}{\updefault}{ $\mu \cdot (n+1)$}%
}}}}
\end{picture}%
\label{picture:construction}

\vspace{0.5cm}


Corollary \ref{continuous extension} tells us that the construction of our array of models
as an increasing sequence of towers is
possible in successor cases. In the limit case, let $I_\omega =
\bigcup_{m<\omega}I_m$, and simply define 
$(\bar M,\bar a,\bar N)^\omega\in\K^*_{\mu,I_\omega}$ to be the union of the towers $(\bar M,\bar a,\bar N)^n$.

 To see that the construction satisfies our requirements, 
first notice that the last column of the array,
$\langle M^n_{i_\theta}\mid n<\omega\rangle$, witnesses that
$M^*$ is a $(\mu,\omega)$-limit model. In light of Theorem~\ref{relatively
full is limit} we need only verify that the last row of the array is a
relatively full tower of cofinality
$\theta$.

\begin{claim}
$(\bar M,\bar a,\bar N)^\omega$ is full relative to
$(M^n_i)_{n<\omega,i\in I_\omega}$.
\end{claim}
\begin{proof}
Given $i$ with $i_\alpha\leq i< i_{\alpha+1}$,
let $(p,M^n_i)$ be some strong type in $\St(M^\omega_i)$.  Notice that
by monotonicity of non-splitting
$(p\restriction M^{n+1}_i,M^n_i)\in\St(M^{n+1}_i)$.  By construction
there is a $j\in I_{n+1}$ with $i< j<i_{\alpha+1}$ such that
$(\tp(a_{j}/M^{n+2}_{j}),N^{n+2}_{j})$ is parallel to
$p\restriction M^{n+1}_i$.  We will show that
$(\tp(a_{j}/M^\omega_{j}),N^\omega_{j})$ is parallel to $(p,N)$.

First notice that $\tp(a_{j}/M^\omega_{j})$ does not $\mu$-split over
$N^\omega_{j}=N^{n+2}_{j}$ because $(\bar M,\bar a,\bar N)^\omega$ is a tower.  
Since
$(\tp(a_{j}/M^{n+2}_{j}),N^{n+2}_{j})$ is parallel to
$(p\restriction M^{n+1}_i,M^n_i)$ there is $q\in\gaS(M^\omega_{j})$
such that $q$ extends both $p\restriction M^{n+1}_i$ and
$\tp(a_{j}/M^{n+2}_{j})$.  
By two separate applications of the
uniqueness of non-$\mu$-splitting extensions we know that $q\restriction
M^\omega_i=p$ and
$q=\tp(a_{j}/M^{\omega}_{j})$.  
To see that $(q,N^{\omega}_{j})$ is parallel to $(p,M^n_i)$, let $M'$ be
an extension of $M^\omega_{j}$ of cardinality $\mu$.  Since
$(p\restriction M^{n+1}_i,M^n_i)$ and $(q\restriction
M^{n+2}_{j},N^{n+2}_{j})$ are parallel, there is $q'\in\gaS(M')$
extending both $p\restriction M^{n+1}_i$ and $q\restriction
M^{n+2}_{j}$ and not $\mu$-splitting over both $M^n_i$ and
$N^{n+2}_{j}$.  By the uniqueness of non-$\mu$-splitting extensions, we
have that $q'$ is also an extension of $q$ and $p$.  Thus $q'$ witnesses
that $(q,N^{\omega}_{j})$ and $(p,M^n_i)$ are parallel.
\end{proof}

This completes the proof of Theorem \ref{main Theorem}.

\section{Concluding remarks}\label{context analysis}

In this section we discuss other results related to Question
\ref{q:Uniqueness problem}.
First to understand the boundaries of Question \ref{q:Uniqueness problem}, consider the elementary case.  Limit models are not necessarily unique even for first order complete
stable theories.

\begin{theorem}\label{strictlystablenonunique}
  Suppose $T$ is a complete, stable theory. Let $\mu\geq 2^{|T|}$ such that
  $\mu^{|T|}=\mu$. If $T$ is not superstable, then no $(\mu,\omega)$-limit model is isomorphic to any
 $(\mu,\kappa)$-limit model
 for any $\kappa$ with $\cf(\kappa)\geq \kappa(T)$.
\end{theorem}

\begin{proof}  Let $T$ be a stable, but not superstable, complete theory, and fix $\kappa$ and $\mu$ as in the statement of the theorem.
  As $T$ is  not superstable, by~\cite[Lemma VII, 3.5 (2)]{Sh e}, for
  $\lambda:=(2^\mu)^+$, there are $\langle \bar{a}_\eta|\eta\in{}^{\omega\, \, \geq}\lambda\rangle$ and $\langle
  \varphi_n(\bar{x},\bar{y}_n)|n<\omega\rangle$ such that
 for every $n<\omega$ and all $\eta\in{}^\omega\lambda$, 
 $$(\C\models \varphi_n[\bar{a}_\eta,\bar{a}_\nu])\Longleftrightarrow
  \nu=\eta\restriction n.$$

  By induction on $n<\omega$ define $\langle M_n|n<\omega\rangle$ all of cardinality $\mu$ and
  $\langle \eta_n,\nu_n|n<\omega\rangle$ such that
  \begin{enumerate}
  \item $M_{n+1}$ is universal over $M_n$ and saturated of cardinality
    $\mu$,
  \item $\eta_{n+1} > \eta_n$, $\nu_{n+1}>\eta_n$,  and $\eta_{n+1}\not= \nu_{n+1}$,
  \item $\bar{a}_{\eta_{n+1}},\bar{a}_{\nu_{n+1}}\in M_{n+1}$ and
  \item \label{equal types condition} $\ftp(\bar{a}_{\eta_{n+1}}/M_n)=\ftp(\bar{a}_{\nu_{n+1}}/M_n)$.
  \end{enumerate}

{\sc This construction is enough:}  Let $N'\models T$ be a $(\mu,\kappa)$-limit over $M_0$.
  By Theorem~\ref{limitfirstorder}, $N'$ must be saturated.  Let $N=\bigcup_{n<\omega}M_n$. Clearly $N$ is a 
  $(\mu,\omega)$-limit over $M_0$.  To conclude that $N$ and $N'$ are non-isomorphic, it is
  enough to show that $N$ is not saturated. Consider $p:=\{
  \varphi_{n+1}(\bar{x};\bar{a}_{\eta_{n+1}})
  \land\lnot\varphi_{n+1}(\bar{x};\bar{a}_{\nu_{n+1}})|n<\omega\}$. The set of formulas $p$
  is a type since it is realized in $\mathfrak C$ by $\bar{a}_\eta$ where
  $\eta:=\bigcup_{n<\omega}\eta_n$. 
  Notice that $N$ cannot satisfy $p$.  If $\bar a\in N$ would satisfy $p$, then $M_n$ realizes $p$ for some $n<\omega$.
  Thus by condition (\ref{equal types condition}), we would have
  \[ {\mathfrak C}\models \varphi_{n+1}[\bar{a},\bar{a}_{\eta_{n+1}}]\Longleftrightarrow {\mathfrak C}\models
  \varphi_{n+1}[\bar{a},\bar{a}_{\nu_{n+1}}]\]
  which would contradict the assumption that $\bar{a}$ satisfies $p$.
    
  {\sc This is possible:} By stability and $\mu^{|T|}=\mu$,  using
  the proof of~\cite[Th. III 3.12]{Sh e}, every model of cardinality
  $\mu$ has a saturated proper elementary extension. Let $M_0$ be a saturated
  model of cardinality $\mu$ and take $\eta_0=\nu_0:=\langle \rangle$. Given $\eta_n,\nu_n,M_n$, using
  Theorem~\ref{exist univ} let $M^*$ be universal over $M_n$ of
  cardinality $\mu$. Let $M^{**}\succ M^*$ of cardinality $\mu$ containing
  $\bar{a}_{\eta_n}$ and $\bar{a}_{\nu_n}$. By~\cite[Th. III 3.12]{Sh e},
  we can take $M_{n+1}\succ M^{**}$ saturated of cardinality $\mu$. Clearly
  it is universal over $M_n$. For $n<\omega$, consider $F_n(\alpha):={\rm tp}(\bar{a}_{\eta_n\conc
  \alpha}/M_n)$. As $\lambda$ is regular and $\lambda>|S(M_n)|$, there is $S\subset \lambda$ of
  cardinality $\lambda$ such that $\alpha\not= \beta\in S \Rightarrow F_n(\alpha)=F_n(\beta)$. Pick $\alpha\not= \beta\in
  S$ and define $\eta_{n+1}:=\eta_n\conc \alpha$ and $\nu_{n+1}:=\eta_n\conc \beta$.
\end{proof}

In the non-elementary setting, many authors have considered approximations to Theorem
\ref{main Theorem}.  Several authors have proved and used the uniqueness of limit models in AECs under the assumption of categoricity:  \cite{Sh 394} \cite{Ba}, \cite{KoSh}, \cite{Sh576}, \cite{ShVi}, \cite{Va1}, and \cite{Va2}. 
Also, Shelah's \cite{Sh i} examines (as an aside) the uniqueness of limit models in good frames.  Below we briefly describe the results and techniques of these papers and distinguish them from our context.

In Theorem 6.5 of \cite{Sh 394}, Shelah claims uniqueness of limit
models of cardinality $\mu$ for classes with
the amalgamation property under little more than categoricity in some
$\lambda>\mu>\LS(\K)$ together with existence of arbitrarily large
models.  
Shelah's claim in Theorem~6.5
of~\cite{Sh 394} (isomorphism over the base) seems too strong for the
proof that he suggests.  Instead, he proves that $(\mu,\kappa)$-limit models
are Galois saturated, which implies uniqueness only over models of
size $<\mu$.
The argument in \cite{Sh 394} depends in a crucial way on an
analysis of Ehrenfeucht-Mostowski models.  In our paper, we cannot employ 
Ehrenfeucht-Mostowski machinery because we do not assume here categoricity  or the existence of models above
the Hanf number.

Under similar categoricity assumptions as those in \cite{Sh 394}, more recently, Baldwin
in~\cite{Ba} (Chapter 11) has used methods based on~\cite{Sh 394} to
prove that if $M_1$ and $M_2$ are $(\mu,\sigma_1)$- and $(\mu,\sigma_2)$-limit models
over $N$, respectively, then $M_1\cong M_2$.  Baldwin, however, does not prove that $M_1$
and $M_2$ are isomorphic {\em over $N$}. Our result is therefore much
stronger than that in~\cite{Ba}. 

Kolman and Shelah  in \cite{KoSh} prove the uniqueness of limit models
of cardinality $\mu$
in $\lambda$-categorical AECs that are axiomatized by a
$L_{\kappa,\omega}$-sentence where $\lambda>\mu$ and $\kappa$ is a
measurable cardinal.  
Then Kolman and Shelah use this uniqueness result to prove that
amalgamation occurs below 
the categoricity cardinal in $L_{\kappa,\omega}$-theories with
$\kappa$ measurable.
Both the measurability of
$\kappa$ and the categoricity are used integrally in their proof of
uniqueness.

Shelah in \cite{Sh576} (see Claim
7.8) proved a special case of the uniqueness of limit models under the 
assumption of
$\mu$-AP, categoricity in $\mu$ and in $\mu^+$ as well
as assuming $K_{\mu^{++}}\neq\emptyset$.  In that paper Shelah needs to
produce \emph{reduced types} and use some of their special properties.

In \cite{ShVi}, Shelah and Villaveces attempted to prove a uniqueness
theorem without assuming
any form of amalgamation; however, they assumed that $\K$ is categorical in
some sufficiently large $\lambda$, that every model in $\K$ has a proper
extension and that $2^{\lambda}<2^{\lambda^+}$.
VanDieren  in \cite{Va1} and \cite{Va2} managed to prove the
uniqueness statement
under the assumptions of \cite{ShVi} together with the additional
assumptions that the class is categorical in $\mu^+$ and
 $\K^{am}:=\{M\in\K_\mu\mid M\text{
is an amalgamation base}\}$ is closed under unions of increasing
$\prec_{\K}$ chains.

In \cite{Sh i} the most important new concept is that of a
\emph{$\lambda$-good frame},
which is an axiomatization of the notion of superstability, {\em with
  hypothesis on just one cardinal $\lambda$}. Its full definition is more
than a page long. Shelah's assumptions on the AEC include, among other
things, the amalgamation property, the existence of a forking like
dependence relation and of a family of types playing a role akin to
that of regular types in first order superstable theories -- Shelah
calls them \emph{bs}-types. One of the
axioms of a good frame is the existence of a non-maximal
super-limit model.  This axiom along with $\mu$-stability implies the
uniqueness of limit models of cardinality $\mu$. In Lemma II.4.8
of~\cite{Sh i} he states that in a good frame, limit models are
unique.  (While we
don't claim that we understand Shelah's proof or believe in its
correctness, he explicitly uses the interplay between \emph{bs}-types
and the forking notion as well as no long forking chains and
continuity of forking.)

The formal differences between our
approach and Shelah's \cite{Sh i} can be summarized as follows: 
\begin{enumerate}
\item  Suppose that $\K$ is an AEC with no maximal models satisfying
  the disjoint amalgamation
property over limit models and is categorical in
$\lambda^+$ for some $\lambda>\LS(\K)$; we then get uniqueness of limit
models.
By way of comparison, in order to get
a uniqueness of  limit models, Shelah needs results of \cite{Sh576} (a
99 pages-long
paper) and significant parts of his book
\cite{Sh i} along with the stronger assumptions of
 categoricity in several consecutive cardinals together with several
 additional set-theoretic axioms.  All our results are in
ZFC.


\item  When specialized to the case where $\K$ is the class of models of a
complete first
order theory $T$, Shelah's proof in~\cite[Lemma II.4.8]{Sh i} really
uses the full power of assuming that $T$ is {\em super}stable. The
proof of uniqueness in this paper just needs, in addition to the
stability and unidimensionality of $T$, no splitting chains of length
$\omega$. As the main
interest of our theorem is for the general case of AEC, rather than
just for first order theories, the difference between this paper
and~\cite[Lemma II.4.8]{Sh i} is clearer when understood in light of
the greater picture.

\end{enumerate}


We are particularly interested in Theorem \ref{main Theorem} not only for
the sake of generalizing Shelah's result from \cite{Sh576} but due to the
fact that the first and second author originally used an earlier draft of this uniqueness theorem (which did not assume unidimensionality) along with tools from \cite{Sh 394} in a
crucial step to prove:

\begin{theorem}[Upward categoricity theorem, \cite{GrVa2}\footnote{Some time after Grossberg and VanDieren announced Theorem~\ref{up cat
  thm}, Baldwin circulated an alternative proof of Theorem~\ref{up cat
  thm} that eventually appeared in~\cite{Ba}. Lessmann in~\cite{Les05}
proved the result for $\K$ with $\LS(\K)=\aleph_0$ beginning with categoricity
in $\aleph_1$.}]\label{up cat thm}
Suppose that $\K$ has arbitrarily large models, is $\chi$-tame and
satisfies the amalgamation and joint embedding properties. Let $\lambda$
be such that
$\lambda>\LS(\K)$ and $\lambda\geq\chi$. If 
$\K$ is categorical in 
$\lambda^+$ then
$\K$ is categorical in all $\mu\geq\lambda^{+}$.
\end{theorem}

After the addition of the unidimensionality assumption in 2014 to
resolve an error found in 2012 in the proof of Theorem \ref{reduced
  are cont}, Grossberg and VanDieren have revisited the proof of
Theorem \ref{up cat thm} to insure that the upward categoricity
transfer still holds \cite{GrVa3}.  Grossberg and VanDieren's initial
use of the uniqueness of limit models in
this theorem hints at a connection between classical definitions of
superstability in first order logic and the uniqueness of limit
models.  This link is explored in further work of VanDieren \cite{Va3}.

It is worth mentioning that the links between classical notions of
superstability from first order logic and the uniqueness of limit
models have also produced interesting insights in the connections
between ``continuous model theory'' and so-called ``metric AECs''. The
work of Villaveces and Zambrano~\cite{ViZa} has extended notions of
independence akin to those used here to the metric AEC context, and
at the same time explored various consequences of assuming forms of
uniqueness of limit models in that metric (continuous)
context.


\begin{thebibliography}{10}

\bibitem[Ba]{Ba} Baldwin, John T.  {\bf Categoricity}. University
  Lecture Series. American
  Mathematical Society 50 (2009).    


\bibitem[Dr]{Dr} Drueck, Fred.   Limit Models, Superlimit Models, and Two Cardinal Problems in
Abstract Elementary Classes.  Ph.D. Thesis at the University of Illinois at Chicago.  (2013).

\bibitem[Gr1]{Gr1} Grossberg, Rami. {\bf A Course in Model Theory I},   to be published by Cambridge University Press.


\bibitem[Gr2]{Gr2} 
Grossberg, Rami.
``Classification Theory for Non-elementary Classes."
{\em Contemporary Mathematics}, {\bf 302} (2002): 165--204.   

\bibitem[GrLe] {GrLe} Grossberg, Rami and Olivier Lessmann.
``Shelah's Stability Spectrum and Homogeneity Spectrum."
     {\em  Archive for
mathematical Logic}, {\bf 41}.1 (2002): 1--31.   

\bibitem[GrVa0]{GrVa0}
Grossberg, Rami and Monica VanDieren.
\newblock  ``Shelah's Categoricity Conjecture from a
 Successor for Tame Abstract Elementary Classes." \emph{Journal of
Symbolic Logic}.
{\bf 71}.2 (2006):  553--568.   

\bibitem[GrVa1]{GrVa1}
---.
``Galois-stability for Tame Abstract Elementary Classes."  \emph{Journal of Mathematical Logic}  {\bf  6}.1 (2006): 25--49.   

\bibitem[GrVa2]{GrVa2} ---.
``Categoricity from one successor cardinal in Tame Abstract Elementary Classes."  {
 \emph{Journal of Mathematical Logic} \bf 6}.2 (2006): 181-201.   

\bibitem[GrVa3]{GrVa3} ---.
``Addendum to `Categoricity from one successor cardinal in Tame Abstract Elementary Classes. {
 \emph{Journal of Mathematical Logic} \bf 6}.2 (2006): 181-201.''   Preprint.

 
 \bibitem[HyKe]{HyKe}
 Hyttinen, Tapani and Meeri Kes\"al\"a.
 \newblock ``Independence in Finitary Abstract Elementary Classes."
 \emph{Annals of Pure and Applied Logic}  {\bf 143} (2006):  103--138.   
 
\bibitem[Jo]{Jo}
 J\'{o}nsson, Bjarni.
\newblock ``Homogeneous universal systems."
\newblock {\em Math. Scand.} {\bf 8} (1960):  137--142.   


\bibitem[KoSh]{KoSh}
 Kolman, Oren and Saharon Shelah.
\newblock ``Categoricity of Theories in $L_{\kappa, \omega}$ when
$\kappa$ is a Measurable Cardinal.  Part I."
\newblock {\em Fundamentae Mathematicae} {\bf151} (1996): 209--240.   

\bibitem[Les05]{Les05}
  Lessmann, Olivier.
  \newblock ``Upward Categoricity from a Successor Cardinal for an
  Abstract Elementary Class with Amalgamation."
\newblock {\em Journal of Symbolic Logic} {\bf 70} (2005): 639--661.   

\bibitem[Po]{Po}
Poizat, Bruno. 
\newblock {\bf A course in model theory}.  Springer.  (2000).   

\bibitem[Sh e]{Sh e} Shelah, Saharon. {\bf Classification Theory}. Revised
  edition. North Holland  (1990).   

\bibitem[Sh i]{Sh i} ---. {\bf Classification Theory for
  Abstract Elementary Classes}. (Studies in Logic: Mathematical Logic
  and Foundations). College Publications (2009).    

\bibitem[Sh 3]{Sh 3} ---.
\newblock ``Finite Diagrams Stable in Power."
\newblock {\em Ann. Math. Logic} {\bf 2},(1970): 69--118.   


\bibitem[Sh 48]{Sh 48}
Saharon Shelah.
\newblock Categoricity in $\aleph _{1}$ of sentences 
in $L_{\omega _{1},\omega}(Q)$, \emph{Israel J Math} {\bf 20} (1975):
127-148. 


\bibitem[Sh 87b]{Sh87b} ---.
\newblock {``Classification Theory for Nonelementary Classes. I. The Number of
  Uncountable Models of $\psi \in L_{\omega _{1},\omega }$. Part B},
\newblock {\em {Israel Journal of Mathematics}} {\bf 46} (1983):  241--273.   

\bibitem[Sh 394]{Sh 394}
---.
\newblock ``Categoricity of Abstract Classes with Amalgamation."
\newblock \emph{Annals of Pure and Applied Logic} {\bf 98} (1999): 261--294.   

\bibitem[Sh 472]{Sh 472} ---.  ``Categoricity of Theories in $L_{\kappa^*,\omega}$ when $\kappa^*$ is a Measurable Cardinal.  Part II."  Dedicated to the memory of Jerzy Los. 
{\em Fundamenta Mathematica}  {\bf 170} (2001): 165--196.   

\bibitem[Sh 576]{Sh576}  ---.
\relax  {``Categoricity of an Abstract Elementary Class in Two Successive
  Cardinals."}
\relax {\em {Israel Journal of Mathematics}} {\bf 126} (2001): 29--128.   


\bibitem[Sh 600]{Sh 600}
---.
\newblock Categoricity in abstract elementary classes: going up inductive
step.
\newblock Preprint.   arXiv:math.LO/0011215 (November 2000 version), 82
pages.


\bibitem[Sh 705]{Sh 705}
---.
\newblock Toward Classification Theory of Good $\lambda$ Frames and
Abstract Elementary Classes.   Preprint.


\bibitem[ShVi]{ShVi}
Shelah, Saharon and Andr\'{e}s Villaveces.
\newblock ``Categoricity in abstract elementary classes with no maximal
models."
\newblock \emph{Annals of Pure and Applied Logic}, {\bf 97}.1-3 (1999):
1--25.  



\bibitem[Va1]{Va1}
VanDieren, Monica.
\newblock ``Categoricity in Abstract Elementary Classes with No Maximal
Models."
\emph{Annals of Pure and Applied Logic} {\bf 141} (2006): 108--147.  Print.


\bibitem[Va2]{Va2}
VanDieren, Monica.  \newblock``Erratum to `Categoricity in abstract elementary classes with no maximal models' [Ann. Pure Appl. Logic 141 (2006) 108--147.0''  \emph{Annals of Pure and Applied Logic}
{\bf 164} (2013):  131-133.  Print. 


\bibitem[Va3]{Va3}
VanDieren, Monica.  \newblock``Superstability and Limit Models in Tame
Abstract Elementary Classes.''  Preprint.

\bibitem[ViZa]{ViZa}
Villaveces, Andr\'es and Zambrano, Pedro. \newblock ``Around
independence and domination in metric abstract elementary classes:
assuming uniqueness of limit models.''
\emph{Mathematical Logic Quarterly} {\bf 60} (2014): 211-227. 


\bibitem[Za]{Za}
Zambrano, Pedro.
\newblock Cats, docilidad, y la propiedad de amalgamaci\'on disyunta.
\newblock Preprint.


\bibitem[Zi]{Zi}
Zilber, Boris.  
\newblock ``A categoricity theorem for quasi-minimal excellent classes.''
\newblock \emph{Contemporary Mathematics} {\bf 380} (2005):297--307.
 
\end{thebibliography}
\end{document}